\documentstyle[]{article}
\def\hybrid{\topmargin 0pt      \oddsidemargin 0pt
        \headheight 0pt \headsep 0pt
        \textwidth 160true mm       
        \textheight 231true mm         
        \marginparwidth 0.0in
        \parskip 0pt plus 1pt   \jot = 1.5ex}
\catcode`\@=11
\def\marginnote#1{}

\newcount\hour
\newcount\minute
\newtoks\amorpm
\hour=\time\divide\hour by60 \minute=\time{\multiply\hour by60
\global\advance\minute by-\hour}
\edef\standardtime{{\ifnum\hour<12 \global\amorpm={am}%
        \else\global\amorpm={pm}\advance\hour by-12 \fi
        \ifnum\hour=0 \hour=12 \fi
        \number\hour:\ifnum\minute<10 0\fi\number\minute\the\amorpm}}
\edef\militarytime{\number\hour:\ifnum\minute<10
0\fi\number\minute}

\def\draftlabel#1{{\@bsphack\if@filesw {\let\thepage\relax
   \xdef\@gtempa{\write\@auxout{\string
      \newlabel{#1}{{\@currentlabel}{\thepage}}}}}\@gtempa
   \if@nobreak \ifvmode\nobreak\fi\fi\fi\@esphack}
        \gdef\@eqnlabel{#1}}
\def\@eqnlabel{}
\def\@vacuum{}
\def\draftmarginnote#1{\marginpar{\raggedright\scriptsize\tt#1}}

\def\draft{\oddsidemargin -.5truein
        \def\@oddfoot{\sl preliminary draft \hfil
        \rm\thepage\hfil\sl\today\quad\militarytime}
        \let\@evenfoot\@oddfoot \overfullrule 3pt
        \let\label=\draftlabel
        \let\marginnote=\draftmarginnote
   \def\@eqnnum{(\theequation)\rlap{\kern\marginparsep\tt\@eqnlabel}%
\global\let\@eqnlabel\@vacuum}  }
\relax

\hybrid

\title{G\'eom\'etrie de l'espace tangent sur l'hyperbolo\"\i de quantique}
\author{Par P.~Akueson\\
I.S.T.V., Universit\'e de Valenciennes, France}

\def\bbbr{{\rm I\!R}}       
\def\bbbc{{\rm I\!\!\! C}}  
\def\bbbn{{\bf  N}}         
\def\bbbk{{\bf K}}          

\def\Fun{\rm Fun\, }

\def\Im{\rm Im\, }
\def\Ker{\rm Ker\, }

\def\Sym{\rm Sym\,}
\def\Vect{\rm Vect\,}
\def\uqs{ U_q(sl(2))}

\def\uqsm{ U_q(sl(2))-Mod}

\def\us{ U(sl(2))}

\def\ovl{sl(2)_q}
\def\uovl{U(sl(2)_q)}
\def\ovk{\!\overline{K}}

\def\ovn{\!\overline{N}}

\def\osf{\!\overline{\sf V}}
\def\ovu{\!\overline{U}}
\def\ovv{\!\overline{V}}
\def\ovw{\!\overline{W}}

\def\ot2{\otimes 2}
\def\ot{\otimes }

\def\op{\oplus}

\def\uq{U_q(\sfg)}
\def\sfg{\sf g}
\def\lbc{[}
\def\rbc{]}

\newtheorem{thm}{Th\'eor\`eme}[section]
\newtheorem{prop}{Proposition}[section]
\newtheorem{lemme}{Lemme}[section]

\newtheorem{defin}{D\'efinition}[section]

\begin{document}

\maketitle
\begin{abstract}
We introduce the tangent space on a quantum hyperboloid. We
define an action of this tangent space on the corresponding
\lq\lq{quantum function space}" $\,{\cal A}\,$, what converts the
elements of the tangent space into \lq\lq{braided vector fields}".
The tangent space is shown to be a projective $\,{\cal A}$-module
and we define a quantum (pseudo)metric and a  quantum connection
(partially defined) on it.
\end{abstract}

\section{introduction}

Il est bien connu depuis les ann\'ees 60 que toutes les
constructions principales de la g\'eom\'etrie classique se
g\'en\'eralisent dans le cadre de la superth\'eorie. Dans les
ann\'ees 80, plusieurs notions et constructions de la
g\'eom\'etrie classique et de la superth\'eorie s'\'etendent au
cas de la g\'eom\'etrie tress\'ee (\lq\lq{braided geometry}"). En
particulier il existe une g\'en\'eralisation naturelle de la
notion d'alg\`ebre commutative (ou super-commutative) li\'ee aux
op\'erateurs de tresses involutifs (appel\'es aussi sym\'etries).
En effet  soit $\,A\,$ une alg\`ebre associative de produit
not\'e $\displaystyle\,\circ\,$ et munie d'un op\'erateur de
tresse\footnote{Nous appelons op\'erateur de tresse (ou un
tressage) une solution de l'\'equation de Yang-Baxter quantique
(EYBQ) : $$S^{12}\,S^{23}\,S^{12}\,= \,S^{23}\,S^{12}\,S^{23}$$
o\`u $\,S^{12}\,=\, S\ot id,\,\,S^{23}\,=\,id\ot S $.}
 $\,S\,$
involutif (i.e. $\,S^2 =id\,$) $$\,S\,:\,A^{\ot 2}\,\to\,A^{\ot
2},\,$$ $\,A\,$ est appel\'ee une {\it alg\`ebre S-commutative}
si :
\begin{enumerate}
\item $\,\,\,\,\displaystyle\,\circ\,=\,\circ\,S\,,\,\,$
\item $\,\,\displaystyle\,S\,{\circ}^{12}\,= \,{\circ}^{23}S^{12}S^{23}$.
\end{enumerate}
Dans ce cas, certains aspects du calcul diff\'erentiel tress\'e
ont \'et\'e d\'evelop\-p\'es dans \cite{GRR}. Plus
particuli\`erement la notion de champ de vecteurs se
g\'en\'eralise (dans l'esprit de la superth\'eorie) de la fa\c
con suivante : une application  lin\'eaire
$\,\displaystyle\,X\,:\,A\to\, A\,$ est appel\'ee un {\it champ
de vecteurs tress\'es} si elle v\'erifie la $\,S$-analogue de la
r\`egle de Leibniz
\begin{equation}
X\,(f\,\circ\,g)\,=\,X\,(f)\,\circ\,g + \circ\,{ev}^{23}\,
S^{12}\,(X\,\ot\,f\,\ot\,g) \label{san}
\end{equation}
o\`u {\em ev}\, est l'application d'\'evaluation
$\,\displaystyle\,X \ot f\,\to\,X(f)$. Nous supposons ici que
$\,S\,$ peut \^etre \'etendu \`a un tressage (not\'e encore
$\,S\,$) \lq\lq{transposant}" les fonctions et les op\'erateurs
et qu'en outre l'application d'\'evaluation commute avec $\,S\,$
dans le sens suivant $$e\,v^{23}\,S^{12}\,S^{23}\,(X\ot f\ot
g)\,=\,S\,e\,v^{12}\,(X\ot f\ot g)\,$$ (il en est de m\^eme en
rempla\c cant la fonction $\,g\,$ par un op\'erateur $\,Y$).

Ainsi, lorsque $\,A\,$ est une alg\`ebre $\,S$-commutative, les
champs de vecteurs sur une telle alg\`ebre peuvent \^etre
d\'efinis par la $\,S$-analogue de la r\`egle de Leibniz
(\ref{san}). \vskip 0.2truecm

Si l'alg\`ebre $\,A\,$ est plut\^ot munie d'une application de
tresse non involutive (c'est juste le cas de l'alg\`ebre que nous
consid\'erons ici) il n'existe aucune d\'efinition  g\'en\'erale
de l'analogue tress\'ee d'une alg\`ebre commutative. (Comme le
montre plusieurs exemples y compris celui  de l'hyperbolo\"\i de
quantique, l'application na\"\i ve des relations de la
$\,S$-commutativit\'e (\c ci-dessus) n'est plus raisonnable).
Autrement dit, il n'est pas \'evident de dire ce qu'est une
alg\`ebre $\,S$-commutative ni de d\'ecrire (dans l'esprit de la
superth\'eorie) les \'el\'ements de la g\'eom\'etrie tress\'ee.
Par exemple la S-analogue de la r\`egle de Leibniz d\'efinie par
la relation (\ref{san}) n'est plus valable et donc, il n'est pas
trivial de d\'ecrire l'espace tangent d'une vari\'et\'e
quantique. (Soulignons que pour les vari\'et\'es quantiques
li\'ees aux applications de tresses involutives, l'espace tangent
peut \^etre d\'ecrit en termes de champs de vecteurs tress\'es.)
\vskip 0.2truecm

Le but de cet article est de d\'efinir l'espace tangent sur
l'hyperbolo\"\i de quantique, d'\'etudier ses propri\'et\'es et
d'introduire les analogues tress\'es de certaines structures
classiques (m\'et\-rique, connexion). Soulignons que nous
traitons (dans l'esprit de \cite{S1}) l'espace tangent comme un
module projectif de rang fini sur l'alg\`ebre des
\lq\lq{fonctions polyn\^omes sur l'hyperbolo\"\i de quantique}".
Rappelons que dans \cite{S1}, Serre a \'etablit une correspondance
biunivoque entre les fibr\'es vectoriels et les modules
projectifs de rang fini sur l'anneau des  coordonn\'ees.(Une
analogue de cette correspondance pour les vari\'et\'es lisses
compactes a \'et\'e \'etablie par R.G.Swan.) \vskip 0.2truecm

L'hyperbolo\"\i de quantique (ou tress\'e) est le plus simple
exemple d'{\it orbite quantique}. Par orbite quantique, nous
entendons les alg\`ebres qui v\'erifient les propri\'et\'es
suivantes :
\begin{enumerate}
\item elles sont $\,\uq$-covariantes  (i.e. leur produit $\,\displaystyle\,\circ\,$ est
 $\,\uq$-cova\-riant
\begin{equation}
X\,\circ\,(a\ot b)\,=\,\circ\,\Delta\,X\,(a\ot b),\qquad
\forall\,\,X\,\in\,\uq \label{uqv}
\end{equation}
o\`u a et b sont des \'el\'ements de l'alg\`ebre consid\'er\'ee),
\item elles repr\'esentent une d\'eformation plate
\footnote{Rappelons que l'alg\`ebre $\,\displaystyle\,{\cal
A}_{\hbar}\,$ (o\`u $\,\hbar\,$ est un param\`etre formel) est une
d\'eformation plate de  l'alg\`ebre  $\,\displaystyle\,{\cal
A}_0,$ si on a $${\cal A}_0={\cal A}_\hbar/\hbar{\cal A}_{\hbar}$$
et $\,\displaystyle\,{\cal A}_0 \ot {\bbbk}\lbc \lbc \hbar\rbc
\rbc\,$ est isomorphe \`a $\,\displaystyle\,{\cal A}_{\hbar}\,$
comme  des $\,\displaystyle\,{\bbbk}\lbc \lbc \hbar\rbc
\rbc$-modules (ici le produit tensoriel est compl\'et\'e en
topologie $\,\hbar$-adique).} de leurs analogues classiques :
i.e. des orbites habituelles dans $\,{\sfg}^*\,$ (plus
pr\'ecis\'e\-ment, les alg\`ebres de fonctions sur de telles
orbites),
\item elles sont en un certain sens commutatives.
\end{enumerate}
 La derni\`ere propri\'et\'e est la plus difficile \`a traiter,
car les orbites quantiques font partie d'une famille d'alg\`ebres
$\,\uq$-covariantes et ce n'est pas toujours facile de distinguer
dans cette famille une alg\`ebre \lq\lq{commutative}" au sens
tress\'e  (voir \cite{DGK} o\`u ce probl\`eme est discut\'e).
Pour expliquer le probl\`eme d\'ecrivons d'abord les objets
infinit\'esimaux sur de telles alg\`ebres. \vskip 0.2truecm

Soit $\,{\cal M}\,$ une vari\'et\'e lisse (en particulier une
orbite dans $\,{\sfg}^*\,$) et $$\rho\,:\,\sfg\,\to\,\Vect ({\cal
M})\,$$ la repr\'esentation de $\,\sfg\,$ sur l'espace des champs
de vecteurs sur $\,{\cal M}$. Associons \`a la $\,R$-matrice
$$R\,=\,\sum_{\alpha\,\in\,\Delta^+}\,X_{\alpha}\,\wedge\,X_{-\alpha}
\quad
 \in {\wedge}^2 (\sfg)$$
o\`u $\,\displaystyle\,\{
H_{\alpha}\,,\,X_{\alpha}\,,\,X_{-\alpha} \}\,$ est une base de
Cartan-Weyl en normalisation de Chevalley de l'alg\`ebre de Lie
$\,\sfg\,$ et  $\,\displaystyle\, \Delta^+\,$ d\'esigne son
syst\`eme de racines positives (en supposant une d\'ecomposition
triangulaire fix\'ee de $\,\sfg$),
 le crochet de $\,R$-matrice
$$\{f\,,\,g\}_R\,=\,\mu\,<{\rho}^{\ot 2}\,(R)\,,\,df\,\ot\,dg>,
\quad f,\,g\,\in\,\Fun({\cal M})\,$$ o\`u $\,\mu\,$ est le
produit dans $\,\displaystyle\,\Fun({\cal M})$.

Dans le cas g\'en\'eral, ce crochet ne satisfait pas \`a la
relation de Jacobi. Par cons\'equent, il n'est pas un crochet de
Poisson. Mais il l'est sur certaines vari\'et\'es dites de type
$\,R$-matrices dans \cite{GP}. C'est juste le cas des orbites
dans $\,\sfg^*\,$ qui admettent comme d\'eformations plates, des
alg\`ebres   $\,\uq$-covariantes. (Mais il existe quand m\^eme
des orbites dans $\,\sfg^*\,$ qui sont des d\'eformations plates
$\,\uq$-covariantes et sur lesquelles le crochet de Poisson
correspondant est un peu diff\'erent de
$\,\displaystyle\,\{\,,\,\}_R\,$, voir \cite{DGS}.) ll est facile
de voir que sur de telles orbites, le  crochet de
Kirillov-Kostant-Souriau (KKS) not\'e $\,\displaystyle\,
\{\,,\,\}_{KKS}\,$ (qui est la restriction sur $\,\sfg^*\,$ tout
entier du crochet lin\'eaire dit de Lie-Poisson) et le crochet de
$\,R$-matrice  sont compatibles et donc sur elles, ces deux
crochets engendrent ce qu'on appelle le pinceau de Poisson. Ainsi
la famille de crochets d\'efinie par
\begin{equation}
\displaystyle\,\{\,,\,\}_{a\,,\,b}\,=\,a\,\{\,,\,\}_{KKS} +
b\,\{\,,\,\}_R,\label{pdp}
\end{equation}
est une famille de crochets de Poisson.

Si l'on peut dire qu'un crochet de Poisson d\'efinit \lq\lq{une
direction de d\'eformation}", un pinceau de Poisson nous donne un
plan de ces directions. La quantification double (i.e.
simultan\'ee) de ce pinceau de Poisson conduit aux alg\`ebres non
commutatives tress\'ees. Ces alg\`ebres d\'ependent de deux
param\`etres. Leur produit  est toujours $\,\uq$-cova\-riant.
\vskip 0.2truecm

Dans cet article,  un cas particulier de telles alg\`ebres est
\'etudi\'e. Notamment celle qui d\'ecoule de l'hyperbolo\"\i de
plong\'e comme  une orbite (par rapport \`a l'action coadjointe
du groupe $\,SL(2)$) d'un \'el\'ement semi-simple de
$\,\displaystyle\,sl(2)^*\,$ (nous regardons plus
pr\'ecis\'ement  une famille d'orbites). Etant l'orbite d'un
\'el\'ement semi-simple la vari\'et\'e initiale est une
vari\'et\'e alg\'ebrique affine ferm\'ee. Notons que dans ce cas
particulier il n'est pas difficile de distinguer dans la famille
des alg\`ebres $\,\uqs$-covariantes, une alg\`ebre
\lq\lq{commutative}" tress\'ee. \vskip 0.2truecm

Soulignons que l'hyperbolo\"\i de quantique a \'et\'e introduit
par Podl\`es \cite{P} sous le nom de la sph\`ere quantique. En
fait, la sph\`ere quantique est l'hyperbolo\"\i de quantique
munie d'une involution (nous ne regardons pas ici le probl\`eme
d'une d\'efinition raisonnable d'une involution sur
l'hyperbolo\"\i de quantique). Nous nous contentons de
consid\'erer simplement l'hyperbolo\"\i de quantique sans
involution. Toutefois la forme compacte ou non  de la vari\'et\'e
initiale n'a aucune importance, puisque   nous  ne regardons que
les fonctions polyn\^omes restreintes sur cette vari\'et\'e. (Par
un changement de variables convenable nous pouvons passer de
l'alg\`ebre des fonctions sur la sph\`ere \`a celle des fonctions
sur l'hyperbolo\"\i de correspondant et r\'eciproquement.) \vskip
0.2truecm

Certains aspects de l'alg\`ebre provenant de la quantification du
pin\-ceau de Poisson  (\ref{pdp}) sur l'hyperbolo\"\i de  ont
\'et\'e examin\'es dans \cite{DG}, \cite{DGR}, \cite{GV}.  En
particulier il a \'et\'e montr\'e que cette alg\`ebre est sans
multiplicit\'e (i.e. chaque module irr\'eductible de dimension
finie entre dans la d\'ecomposi\-tion en $\,\uqs$-modules
irr\'eductibles de cette alg\`ebre au maximum une fois). \vskip
0.2truecm

D\'esignons  par $\,\displaystyle\,{\cal A}_{\hbar,q}^c\,$ cette
alg\`ebre  \`a deux  param\`etres qui est le r\'esultat de la
quantification double du pinceau de Poisson sur l'hyperbo\-lo\"\i
de plong\'e comme une  orbite dans $\,sl(2)^*$. En gros, nous
disons que le param\`etre $\,\hbar\,$ est celui de la
quantification du crochet de KKS,  $\,q\,$ celui de tressage et
$\,c\,$ num\'erote les orbites quantiques. Dans le cas o\`u
$\,\hbar=0\,$, nous obtenons l'alg\`ebre $\,\displaystyle\,{\cal
A}_{0,q}^c\,$ qui est consid\'er\'ee comme la $q$-analogue de
l'alg\`ebre commutative correspondante. (Elle est parfois
appel\'ee sph\`ere quantique standard.) C'est notre alg\`ebre
principale : \lq\lq{l'alg\`ebre des fonctions polyn\^omes
restreintes sur l'hyperbolo\"\i de quantique}" (o\`u
\lq\lq{l'anneau des coordonn\'ees quantiques}"). Ainsi dans la
suite, tous les modules consid\'er\'es (dans le cas quantique)
seront des modules sur cette alg\`ebre. Nous introduisons alors
les notions de champs de vecteurs, d'espace tangent sur
l'hyperbolo\"\i de quantique. Nous traitons cet espace tangent
comme un $\,{\cal A}_{0,q}^c$-module et nous l'appellons
\lq\lq{module tangent}" sur l'hyperbolo\"\i de quantique. \vskip
0.2truecm

Commen\c cons par le module tangent. Dans le cas classique, il
est form\'e par des champs de vecteurs tangents \`a la
vari\'et\'e sous consid\'eration (ici la sph\`ere ou
l'hyperbolo\"\i de). Quels sont les analogues quantiques de ces
champs ?

Les champs de vecteurs auxquels l'on peut spontan\'ement penser,
sont les g\'en\'erateurs $\,\displaystyle\,X,Y,H\,$ du
 groupe quantique (GQ) $\,\uqs\,$ (voir Section 2).
Mais si nous introduisons le module tangent comme toutes les
combinaisons lin\'eaires (ayant pour coefficients les
\'el\'ements de l'alg\`ebre $\,\displaystyle\,{\cal A}_{0,q}^c$)
 de ces op\'erateurs, nous n'obtenons pas une d\'eformation
plate du module tangent initial. En effet, notons  $\,\Fun(S^2)\,$
l'alg\`ebre des fonctions  polyn\^omes restreintes sur la
sph\`ere ($\,S^2\,$). On peut d\'ecrire le \lq\lq{module
tangent}"  sur la sph\`ere comme toutes les combinaisons
lin\'eaires aux coefficients-fonctions de trois rotations
infinit\'esimales
\begin{equation}
X=z\partial_y - y\partial_z,\,\, Y=x\partial_z - z\partial_x,\,\,
Z=y\partial_x -x\partial_y. \label{eqc}
\end{equation}
Ces rotations  correspondent aux g\'en\'erateurs standards
$\,x,\,y,\,z\,$ de $\,{\sfg}=so(3)=su(2)\,$ qui op\`erent sur
l'alg\`ebre de Lie $\,\sfg\,$ elle-m\^eme par l'action adjointe,
et leur extension sur les \'el\'ements de
$$\,\Fun({\sfg^*})=Sym\,\,({\sfg})\,$$ s'effectue par la r\`egle
de Leibniz. Il est facile de voir que sur l'alg\`ebre
$\,\Fun(S^2)$, les op\'erateurs $\,\displaystyle\,X,\,Y,\,Z\,$
satisfont \`a la relation $$x\,X + y\,Y + z\,Z=0.$$

En passant \`a la forme non compacte (i.e. \`a l'hyperbolo\"\i de
$\,{\rm H}\,$), nous avons la relation
\begin{equation}
\displaystyle\,\,xY + yX + \frac{h}{2}H\,=\,0 \label{yoe}
\end{equation}
avec les g\'en\'erateurs standards $\,x,\,y,\, h\,$ de l'alg\`ebre
$\,sl(2)\,$ et $\,X,\,Y,\, H\,$ sont ici les rotations
hyperboliques infinit\'esimales correspondantes. Notons
$\,\Fun({\rm H})\,$ l'alg\`ebre des fonctions polyn\^omes
restreintes sur l'hyperbo\-lo\"\i de. Finalement nous pouvons
d\'efinir l'espace tangent sur l'hyperbo\-lo\"\i de comme un
module facteur  d'un module libre $\,A^3\,$ de rang  3 (ici
$\,\displaystyle\,A={\cal A}_{0,1}^c=\Fun({\rm H})$) de la fa\c
con suivante : $$\,aX + bY + cH\quad{\rm modulo}\quad
f.(xY+yX+\frac{h}{2}H)\,$$ o\`u $\,a,\,b,\,c,\,f\,\in \,\Fun({\rm
H})$. Notons   $\,T({\rm H})\,$ le $\,\Fun({\rm H})$-module
tangent sur l'hyperbolo\"\i de. Nous avons en outre l'action
$$T({\rm H})\ot A\,\to \,A$$ qui signifie que les \'el\'ements du
$\,A$-module (disons gauche) $\,T({\rm H})\,$ sont pr\'esent\'es
comme des op\'erateurs  sur $\,A\,$ et on appelle le plongement
\begin{equation}
sl\,(2)\quad \hookrightarrow \quad T({\rm H}).\label{plon}
\end{equation}
une {\it ancre}. \vskip 0.2truecm

Notre but est de d\'efinir le module tangent sur l'hyperbolo\"\i
de quantique muni d'une  {\it ancre quantique}. \vskip 0.2truecm

Contrairement au cas classique, les g\'en\'erateurs
$\,X,\,Y,\,H\,$ du GQ $\,\uqs\,$ ne satisfont \`a aucune relation
du type (\ref{yoe}) et donc la platitude de la d\'eformation du
module tangent (engendr\'e par les op\'erateurs $\,X,\,Y,\,H\,$)
n'a pas lieu.

Nous sugg\'erons d'autres candidats pour le r\^ole des
$q$-analogues des op\'erateurs $\,X,\,Y,\,H\,(\in \us)\,$,  de
sorte qu'ils satisfassent \`a la $q$-analogue de la  relation
(\ref{yoe}). Donc l'analogue quantique du module tangent peut
\^etre introduit  de fa\c con analogue au cas classique. La
construction  se fait en deux \'etapes :

$\bullet\,\,$
 Primo, nous d\'efinissons   l'analogue tress\'e du crochet
de Lie de $\,sl(2)\,$ (il a \'et\'e introduit pour la premi\`ere
fois dans \cite{DG} et g\'en\'eralis\'e apr\`es dans \cite{LS}).
Cela nous permet d'introduire l'ad-action de l'alg\`ebre
\lq\lq{$\,sl(2)\,$ {\it tress\'ee}}" sur elle-m\^eme. Il est
facile de v\'erifier que les op\'erateurs correspondants aux
g\'en\'erateurs de l'alg\`ebre $\,sl(2)\,$ tress\'ee,  satisfont
\`a une relation qui est la $q$-analogue de la relation
(\ref{yoe}). Cela incite alors \`a d\'efinir  l'analogue
quantique du module tangent  comme un module quotient d'un module
libre en utilisant  cette relation.

$\bullet\,\,$ Secundo, si nous voulons avoir une action de ce
module sur l'alg\`ebre $\,{\cal A}_{0,q}^c\,$, nous devons
\'etendre l'action des op\'erateurs adjoints qui est bien
d\'efinie pour le moment sur les \'el\'ements de degr\'e un de
l'alg\`ebre $\,{\cal A}_{0,q}^c\,$ \`a toute fonction polyn\^omes
de degr\'e sup\'erieur \`a un, sur l'hyperbolo\"\i de quantique.
\vskip 0.2truecm

Malheureusement pour le cas que nous traitons, il n'existe
aucune forme tress\'ee de la r\`egle de Leibniz.
 Nous pr\'esentons une autre m\'ethode pour \'etendre l'action des
op\'erateurs adjoints provenant de l'ad-action de l'alg\`ebre
$\,sl(2)\,$ tress\'ee, sur toute l'alg\`ebre $\,{\cal
A}_{0,q}^c\,$ de sorte que la $q$-analogue de la relation
(\ref{yoe}) soit toujours satisfaite.

En fait notre construction nous am\`ene \`a  deux modules
tangents : celui trait\'e comme un $\,\displaystyle\,{\cal
A}_{0,q}^c$-module
 gauche not\'e $\,T({\rm H}_q)_l\,$ et celui trait\'e comme un
 $\,\displaystyle\,{\cal A}_{0,q}^c$-module droit not\'e
$\,T({\rm H}_q)_r$. Le probl\`eme est de les identifier,
c'est-\`a-dire construire un isomorphisme (dans la cat\'egorie
consid\'er\'ee\footnote{Par d\'efinition, un  morphisme dans
cette cat\'egorie
 est introduit comme une application
qui commute \`a l'action de $\uqs$.}) entre ces deux modules.
 Notons qu'une telle identification est une \'etape
indispensable pour d\'efinir par exemple une (pseudo)m\'etrique
tress\'ee sur l'hyperbolo\"\i de quantique. Le probl\`eme est que
notre m\'ethode  pour introduire une telle m\'etrique consiste
\`a d\'efinir d'abord un \lq\lq{couplage}"
$\,\displaystyle\,<\,,\,>\,$ sur $\,\displaystyle\,T({\rm H}_q)_l
{\ot}_{\bbbk} T({\rm H}_q)_r\,$ $$<\,,\,>\,:\,T({\rm H}_q)_l
{\ot}_{\bbbk} T({\rm H}_q)_r\,\to \,{\cal A}_{0,q}^c.$$ (Ici,
$\,T({\rm H}_q)_l\,$  est un $\,{\cal A}_{0,q}^c$-module gauche
et  $\,T({\rm H}_q)_r\,$ un $\,{\cal A}_{0,q}^c$-module  droit.)

Le passage \`a une (pseudo)m\'etrique  sur $$T({\rm
H}_q)_r{\ot}_{\bbbk}T({\rm H}_q)_r\,\to\, {\cal A}_{0,q}^c
\qquad  \mbox{ou}\qquad T({\rm H}_q)_l{\ot}_{\bbbk} T({\rm
H}_q)_l\,\to \,  {\cal A}_{0,q}^c$$ ne peut \^etre effectu\'e
qu'apr\`es l'identification des modules $\,T({\rm H}_q)_l\,$ et
$\,T({\rm H}_q)_r$ : $$T({\rm H}_q)_l\approx T({\rm H}_q)_r.$$
Soulignons que sans une telle identification, il n'est pas
\'evident de sortir le facteur $\,f\,$ dans le terme
$$<X\,,\,fY>\quad X\,,Y\,\in \,T({\rm H}_q)_l\,,\,\,
f\,\in\,{\cal A}_{0,q}^c.$$ Plus pr\'ecis\'ement sans cette
identification on ne peut pas d\'efinir le produit tensoriel
$\,T({\rm H}_q)_{\epsilon} {\ot}_{{\cal A}_{0,q}^c} T({\rm
H}_q)_{\epsilon}\,$ o\`u $\,{\epsilon}=l,\,r$.

Nous  construisons \'egalement une connexion (partiellement
d\'efinie, i.e. sur un sous espace de $\,T({\rm H}_q)_{\epsilon}
\ot T({\rm H}_q)_{\epsilon}\,$) tress\'ee. Le probl\`eme dans ce
cas est que, pour \'etendre cette connexion sur
 $\,T({\rm H}_q)_{\epsilon} \ot T({\rm H}_q)_{\epsilon}\,$,
il faut pouvoir \'etendre le crochet de Lie tress\'e sur cet
espace. Malheureusement nous ne savons pas  prolonger (par les
m\'ethodes existantes) sur  $\,T({\rm H}_q)_{\epsilon} \ot T({\rm
H}_q)_{\epsilon}\,$ le crochet de Lie tress\'e initialement
d\'efini  sur $\,sl(2)$. \vskip 0.2truecm

Pr\'ecison  que les crit\`eres de raison d'\^etre pour les objets
et les op\'era\-teurs que nous introduisons sont de deux sortes :
la platitude de la d\'eformation et la $\,\uqs$-covariance.
\vskip 0.2truecm

Le contenu de cet article est le suivant. Dans la Section 2, nous
pr\'esentons notre alg\`ebre principale : celle des fonctions sur
l'hyperbolo\"\i \-de quantique. Nous rappelons dans la Section 3,
l'alg\`ebre  enveloppante tress\'ee de $\,sl(2)$. Dans la Section
4, nous d\'efinissons les {\it champs de vecteurs tress\'es} et
l'espace tangent (consid\'er\'e comme un $\,{\cal
A}_{0,q}^c$-module) sur l'hyperbolo\"\i de quantique muni d'une
ancre quantique. Dans la Section 5 nous montrons que ce $\,{\cal
A}_{0,q}^c$-module tangent
 est projectif. Dans la Section 6 nous
d\'efinissons et construisons une (pseudo)\footnote{Pseudo
signifie que son analogue classique n'est pas d\'efinie
positive.}m\'etrique et une connexion (partiellement d\'efinie)
tress\'ees
 sur l'hyperbolo\"\i de quantique. Dans la derni\`ere
Section, nous proposons une fa\c con canonique pour identifier
les modules
tangents gauche et droit sur l'hyperbolo\"\i de quantique. \\

Dans toute la suite, le corps de base $\,{\bbbk}\,$ est
$\,{\bbbr}\,$ ou $\,{\bbbc}\,$ et le
param\`etre $\,q\,$ ($\in \,\bbbk$) est g\'en\'erique.\\

\vskip 2truecm

{\bf Remerciemments} : D. Gurevich m'a soumis ce probl\`eme et
m'a  constamment guid\'e dans sa r\'esolution et dans sa
r\'edaction. Je l'en remercie.

\newpage

\section{Hyperbolo\"\i de  quantique}

Rappelons d'abord la pr\'esentation de l'hyperbolo\"\i de
classique.

Soient :

\noindent $sl(2)$ l'alg\`ebre de Lie du groupe de Lie  $\,SL(2)$,

\noindent $\,\lbc\,,\,\rbc\,$ le crochet de Lie sur $sl(2)$,

\noindent $\,sl(2)^*\,$ l'espace (vectoriel) dual de $\,sl(2)\,$
et

\noindent $(X,\,Y,\,H)\,$ la base de Cartan-Weyl de $\,sl(2)\,$.

La repr\'esentation \lq\lq{adjointe gauche}" (Ad) de $\,SL(2)$ :
$$Ad_g\,X\,=\,g\,X\,g^{-1},\quad g\,\in\,SL(2),\,\,X\,\in\,sl(2)$$
 fait de $\,sl(2)\,$ un $SL(2)$-module gauche.
(Notons que Ad d\'esigne  souvent la repr\'esentation adjointe
droite. Mais nous pr\'ef\'erons  r\'ealiser $\,sl(2)\,$ comme un
$SL(2)$-module gauche.)

Par la repr\'esentation coadjointe  associ\'ee et not\'ee
$\,Ad^{*}\,$, l'espace $\,sl(2)^*\,$ est un $SL(2)$-module droit.

Soit $\,\displaystyle\,\omega$ l'\'el\'ement de $sl(2)^*$
d\'efini par : $$\,\displaystyle\,\omega (H)=a,\,\, \omega
(X)\,=\,\omega (Y)=0,\,\,a\in \,\bbbk,\,\,a\not=0.$$ D\'esignons
par $\,\displaystyle\,{\cal O}_{\omega}\,$ l'orbite de
l'\'el\'ement $\omega\,$ par  la repr\'esentation coadjointe
$\,Ad^{*}$ : $$\displaystyle\,{\cal
O}_{\omega}\,=\,\{\,{Ad_g}^{*}(\omega)\,\, /
\,\,g\in\,SL(2)\,\}.$$ Le stabilisateur de l'\'el\'ement
$\,\displaystyle\,\omega\,$ est juste le sous-groupe de Cartan
$\,{\bf H}\,$ de $\,SL(2)$. Par cons\'equent comme un espace
homog\`ene :
\begin{eqnarray*}
&\displaystyle\,{\cal O}_{\omega}=SL(2)/{\bf H}.
\end{eqnarray*}

Pr\'esentons maintenant  l'orbite $\,\displaystyle\,{\cal
O}_{\omega}$ comme une  vari\'et\'e alg\'ebrique affine.

Consid\'erons $\,\displaystyle\,\Fun(sl(2)^*)\,$, l'alg\`ebre des
fonctions polyn\^omes sur l'espa\-ce $\,sl(2)^*$. Soit
$\,\displaystyle\,\Sym(sl(2))\,$ l'alg\`ebre sym\'etrique de
l'espace $\,sl(2)$.  On a alors de mani\`ere naturelle :
$$\,\displaystyle\,\Fun(sl(2)^*)=\Sym(sl(2)).$$ Soit $\,\us$
l'alg\`ebre enveloppante de $\,sl(2)$. Elle est une alg\`ebre
filtr\'ee. Soit $\,\displaystyle\,{\rm Gr}\,\us\,$ l'alg\`ebre
gradu\'ee associ\'ee \`a l'alg\`ebre filtr\'ee $\,\us$.  On a :
$\,\displaystyle\,{\rm Gr}\,\us \cong \Sym\,(sl(2))\,$ par le
th\'eor\`eme de Poincar\'e-Birkhoff-Witt (PBW).

Soit $\,{\rm C}\,$ l'\'el\'ement de Casimir de $\,\us\,$
$$\,\displaystyle\,{\rm C}\,=\,XY + YX + \frac{H^2}{2}\,\,,$$
$\rm C\,$, est un g\'en\'erateur du centre de l'alg\`ebre $\us$.
Associons \`a chaque \'el\'ement $\,Z\,$ de $\,\us\,$ son image
(not\'ee $z$) dans $\,\displaystyle\,{\rm Gr}\,\us\,$
($\approx\,\Sym\,(sl(2))$). Alors par cette correspondance,
l'image  (que nous notons encore $\rm C\,$) de l'\'el\'ement de
Casimir  dans $\,\displaystyle\,\Sym(sl(2))\,$ est
$$\,\displaystyle\,{\rm C}\,=\,2\,xy +\,\frac{h^2}{2}.$$

Il est bien connu que toute orbite de la repr\'esentation
coadjointe $\,Ad^{*}\,$ est contenue dans la vari\'et\'e
alg\'ebrique affine d\'efinie par : $$\displaystyle\,{\rm C
}\,=\,2\,xy +\,\frac{h^2}{2}\,=\, c\,\,\,\mbox{o\`u}\quad
c\,\,\,\mbox{est une constante dans} \quad\bbbk.$$ En particulier
si $\,\bbbk\,=\,\bbbc\,$ et $\,c\not=0\,$, l'orbite
$\,\displaystyle\,{\cal O}_{\omega}\,$ co\"\i ncide avec
l'hyperbolo\"\i \-de (classique) $\,{\rm H}\,$ d'\'equation
\begin{equation}
2xy + \frac{h^2}{2}\,=\,c\,=\,{\rm
C}(\omega)\,=\,\frac{a^2}{2}.\label{orc}
\end{equation}
Si  $\,\bbbk\,=\,\bbbr\,$, l'hyperbolo\"\i de contient parfois
une orbite, parfois deux.

Pour $\,c=0\,$, (\ref{orc}) d\'efinit le  c\^one qui est
compos\'e de deux orbites $\,\{0\}\,$ et tout le reste si
$\,\bbbk\,=\,\bbbc\,$; et de trois orbites si
$\,\bbbk\,=\,\bbbr$. Fixons $\,c\not =0\,$ et consid\'erons
$\,\displaystyle\,\Fun({\rm H})\,$, l'alg\`ebre des fonctions
polyn\^omes sur l'hyperbolo\"\i de $\,\displaystyle\,{\rm H}$.
Par d\'efinition  $\,\displaystyle\,\Fun({\rm H})\,$ est la
restriction des fonctions polyn\^omes  de
$\,\displaystyle\,\Fun(sl(2)^*)\,$ sur  $\,\displaystyle\,{\rm
H}\,$ i.e.
\begin{equation}
\Fun({\rm H}) = \Fun(sl(2)^*)/\{{\rm C}-c\},\label{fyo}
\end{equation}
o\`u $\,\displaystyle\,\{{\rm C}-c\}\,$ d\'esigne l'id\'eal
bilat\`ere engendr\'e par  l'\'el\'ement $\,\displaystyle\,{\rm
C}-c\,$.

Il est en outre facile de voir que la multiplication dans
l'alg\`ebre $\,\displaystyle\,\Fun({\rm H})\,$ est  covariante
par l'action de $\,\us$. \vskip 0.2truecm

Par analogie avec le cas classique pr\'ec\'edemment d\'ecrit,
nous pr\'esen\-tons l'analogue  quantique not\'e
$\,\displaystyle\,{\rm H}_q\,$ de l'hyperbolo\"\i de  classique
$\,{\rm H}\,$, sous forme de son alg\`ebre des \lq\lq{fonctions
quantiques}". La multiplication dans cette alg\`ebre doit \^etre
en outre $\,\uqs$-covariante. Pour ce faire, nous  donnons dans
la suite un analogue \lq\lq{tress\'e}" de l'\'el\'ement de
Casimir $\,{\rm C}\,$ qui participera \`a nos constructions.

Soit $\,\uqs\,$ le GQ associ\'e au groupe $\,SL(2)$. Le groupe
 $\,\uqs$
est une alg\`ebre de Hopf. Dans le mod\`ele de Drinfel'd-Jimbo,
elle est engendr\'ee par les \'el\'ements $\,X,\,Y,\,H\,$
satisfaisant aux relations de commutation (pour $q\not =0,\,\,q^2
\not =1$) :
\begin{equation}
\lbc H,X\rbc=2X ,\quad \lbc H,Y \rbc=-2Y ,\quad \lbc X,Y \rbc
={q^H-q^{-H} \over q-q^{-1}} \label{pue}.
\end{equation}
On peut choisir le coproduit ($\,\displaystyle\,\Delta \,$)
d\'efini par exemple par :
\begin{equation}
\Delta (X) = X\otimes 1  +  q^{-H}\otimes  X,\ \ \Delta (Y) =
1\otimes Y  +  Y\otimes q^H,\ \ \Delta (H) = H\otimes 1 +
1\otimes H \label{cuq}.
\end{equation}
Alors l'antipode $\,\displaystyle\,\gamma\,$  est donn\'ee par :
\begin{equation}
\gamma (X)=-q^H X,\,\,\gamma (H)=-H,\,\,\gamma (Y)=-Yq^{-H}.
\end{equation}
\noindent (Pour $\,\displaystyle\,q=1\,$, les relations
(\ref{pue}) et (\ref{cuq}) correspondent \`a celles de
l'alg\`ebre de Hopf $\,\us$). \vskip 0.2truecm

D\'esignons par $\,\uqsm\,$ la cat\'egorie des $\uqs$-modules de
dimension finie qui sont  analogues quantiques (c'est-\`a-dire
des d\'eforma\-tions)  de $\us$-modules irr\'eductibles de
dimension finie. Tout objet de $\,\uqsm\,$ est appel\'e
$q$-analogue (ou analogue tress\'e) de l'objet classique
correspondant. \vskip 0.2truecm

Le centre de l'alg\`ebre $\,\uqs\,$ est engendr\'e (voir
\cite{M}) par l'op\'era\-teur de  {\it Casimir quantique}

\begin{equation}
{\rm
C}_q=(\frac{q^{\frac{H+1}{2}}-q^{-\frac{H+1}{2}}}{q-q^{-1}})^2 +
YX \,.\label{cas}
\end{equation}

\noindent (On peut remarquer que  pour $\,q=1\,$ on a $\,\quad
{\rm C}_1 =\frac{\rm C}{2}+ \frac{id}{4}\not={\rm C}$.)

Donnons \`a pr\'esent un autre analogue de l'\'el\'ement de
Casimir $\,{\rm C}\,$ qui nous servira dans la suite. Pour cela,
consid\'erons  une seconde copie de l'espace $\,sl(2)\,$ que nous
notons $\,{\sf V}\,$ pour la diff\'erencier de l'espace initiale
$\,sl(2)\,$ et d\'esignons par ($\,u,\,v,\,w\,$) une base de
l'espace   $\,{\sf V}$ : $${\sf V}\,=\,Span\,(u,\,v,\,w)\,,\,\,
\mbox{notons}\quad sl(2)\,:=\,({\sf V}\,,\,\lbc\,,\,\rbc).$$

Munissons $\,\displaystyle\,{\sf V}\,$ de l'action de $\,\uqs\,$
qui co\"\i ncide pour $q=1$ avec celle de la repr\'esentation
adjointe (ad) de $\,sl(2)$. Elle est not\'ee {\bf .} et
d\'efinie  par :
\begin{eqnarray*}
\displaystyle\,X.u&=&0,\,\,\,X.v=-(q+q^{-1})u,\,\,\,X.w=v,\\
\displaystyle\,Y.u&=&-v,\,\,\,Y.v=(q+q^{-1})w,\,\,\,Y.w=0,\\
\displaystyle\,H.u&=&2u,\,\,\,H.v=0,\,\,\,H.w=-2w.
\end{eqnarray*}
(Pour $q=1\,$, on v\'erifie qu'on a bien l'\lq\lq{ad-action
gauche}" de $\,sl(2)$.)

La structure de coalg\`ebre de $\,\uqs\,$ permet d'\'etendre
cette action sur $\,\displaystyle\,{\sf V}^{\ot2}$. Par le
th\'eor\`eme de Clebsch-Gordan quantique (voir \cite{K}),
$\,\displaystyle\,{\sf V}^{\ot2}\,\,$ se d\'ecompose en trois
$\,\,\uqs$-modules irr\'eductibles de dimension finie
$\,\displaystyle\,{\sf V}_0^q,\,{\sf V}_1^q,\,{\sf V}_2^q\,$
respectivement de spins $\, 0\,,1\,,2$. Fixons  respectivement
dans les espaces $\,{\sf V},\,{\sf V}_0^q,\,{\sf V}_1^q,\, {\sf
V}_2^q\,$ leurs \'el\'ements de plus haut poids, not\'es :
$$x_0,\quad {\cal C}_q,\quad x_1,\quad x_2$$ o\`u $\,{\cal C}_q =
(q^3 + q) uw + vv + (q + q^{-1})wu$.

La compatibilit\'e de la $\,\uqs$-action sur l'espace $\,{\sf
V}^{\ot 2}\op {\sf V}\op {\bbbk}\,$  impose les relations
\begin{equation}
{\cal C}_q = c,\quad x_1 =\hbar x_0\label{fsc}
\end{equation}
o\`u $\,c\,$ et $\,\hbar\,$ sont des constantes dans $\,{\bbbk}$.
$\displaystyle\,{\cal C}_q\,$ est appel\'e {\it Casimir
tress\'e}. Il n'est pas \`a confondre  avec le Casimir quantique
$\,\displaystyle\,{\rm C}_q\,$ (d\'efini par (\ref{cas})) qui
appartient  \`a l'alg\`ebre  $\,\uqs$.

Faisons remarquer que pour $q=1$, $\,\displaystyle\,{\cal C}_1 =
4uw+v^2=2{\rm C}$. $\,\displaystyle\,{\cal C}_q\,$ (plus
pr\'ecis\'ement $\,\frac{{\cal C}_q}{2}\,$) est  la
$q$-analogue  de $\,\displaystyle\,{\rm C}\,$, dont nous nous
servons dans toute la suite pour nos constructions.

En op\'erant avec l'op\'erateur $\,Y \in \uqs\,$ sur la
deuxi\`eme  \'equation de (\ref{fsc}), on en d\'eduit en plus
deux autres \'equations. (Voir appendice {\bf A} pour leurs
formes explicites). En posant : $$x_0\,=\,u,\quad
x_1\,=\,q^2uv-vu,\quad x_2\,=\,uu$$ on en d\'eduit l'expression
des deux autres \'el\'ements de base de $\,{\sf V}_1^q\,$ en
op\'erant avec $\,Y (\in \uqs)\,$ sur l'\'el\'ement $\,x_1$.

\begin{defin}\label{dey}
L'alg\`ebre $\,\displaystyle\,{\cal A}_{\hbar,q}^c\,$ est le
quotient de l'alg\`ebre  tensorielle libre
$\,\displaystyle\,T({\sf V})\,$ par l'id\'eal bilat\`ere
$\,\displaystyle\,I_{\hbar}\,$ engendr\'e par les \'el\'ements :
\begin{eqnarray*}
&\displaystyle\,q^2uv-vu+2u\hbar,\,\,\,(q^3+q)(uw-wu)+(1-q^2)vv-2v\hbar,\,\,\\
&\displaystyle\,-q^2vw+wv-2w\hbar,\,\,\,{\cal C}_q -c.
\end{eqnarray*}
Pour $\,\hbar=0\quad \mbox{et}\quad c\not=0\,$, l'alg\`ebre
$\,\displaystyle\,{\cal A}_{0,q}^c\,$ est celle des fonctions sur
l'hyperbolo\"\i de quantique (qui est aussi not\'ee $\,{\rm
H}_q$).
\end{defin}
\noindent {\bf Remarque 2.1} Dans cette d\'efinition
$\,\displaystyle\,\hbar\,$ et $\,q\,$ sont des constantes
fix\'ees. Nous pouvons  les consid\'erer comme des param\`etres.
Il suffit pour cela de remplacer dans  la d\'efinition
pr\'ec\'edente $\,\displaystyle\,T({\sf V})\,$ par
$\,\displaystyle\,T({\sf V})\ot{\bbbk}[[\hbar,\,q,\,q^{-1}]]$. Le
param\`etre orbital $\,c\,$ est une constante.

Bien que l'hyperbolo\"\i de classique soit une vari\'et\'e non
compacte, en un certain sens, l'hyperbolo\"\i de quantique est
plut\^ot un analogue tress\'e de la sph\`ere. Ceci vient du fait
que nous ne consid\'erons que les fonctions polynomiales sur
l'objet classique et leurs analogues tress\'es. En outre, toutes
les repr\'esentations que nous consid\'erons sont de dimension
finie.

\begin{prop}
$\,\displaystyle\,{\cal A}_{\hbar,q}^c\,$ est une alg\`ebre
associative $\,\uqs$-covariante.
\end{prop}

$Preuve :$ C'est par construction de l'alg\`ebre
$\,\displaystyle\,{\cal A}_{\hbar,q}^c$. \vskip 0.5truecm

$\bullet$  $\,\,\displaystyle\,{\cal A}_{0,1}^c\,$ (pour
$\,\displaystyle\,c\not=0\,$) est l'alg\`ebre des fonctions
polyn\^omes sur $\,sl(2)\,$ restreintes \`a la vari\'et\'e
alg\'ebrique affine d\'efinie par :
$\,\displaystyle\,4uw+v^2\,=\,c\,$, c'est l'alg\`ebre des
fonctions  sur l'hyperbolo\"\i de (classique) qui est
commutative.

$\bullet$  L'alg\`ebre $\,\displaystyle\,{\cal A}_{\hbar,1}^c\,$
(pour $\,\displaystyle\,c\not=0\,$) est l'analogue non
commutative de l'alg\`ebre $\,\displaystyle\,{\cal A}_{0,1}^c\,$,
mais elle est toujours $\,sl(2)$-invariante.

\vskip 0.5truecm \noindent {\bf Remarque 2.2}\label{rep}
L'alg\`ebre $\,\displaystyle\,{\cal A}_{\hbar,1}^c\,$  est une
d\'eformation plate de  l'alg\`ebre  $\,\displaystyle\,{\cal
A}_{0,1}^c$.
 Cela d\'ecoule du th\'eor\`eme de PBW. Ainsi pour
$\,c\not=0\,$, comme dans la d\'ecomposition en $\,sl(2)$-modules
irr\'eductibles de dimension finie ($\,{\sf V}_k\,$ de spin
$\,k\in {\bbbn}\,$) de l'alg\`ebre $\,{\cal A}_{0,1}^c\,$, toute
composante appara\^\i t  sans multiplicit\'e, il en est de m\^eme
pour les alg\`ebres $\,{\cal A}_{\hbar,1}^c$ et $\,{\cal
A}_{\hbar,q}^c\,$ ($\,q\not=1\,$) (voir \cite{GV},\cite{A}).

$\bullet$  Le cas $\,\displaystyle\,c=0\,$ correspond au c\^one
(dit \lq\lq{quantique}" si $\,\displaystyle\,q\not=1\,$).

$\bullet$ $\,\displaystyle\,{\cal A}_{\hbar,q}^c\,$ est une
famille (d\'ependante de $\,\hbar$) d'alg\`ebres
$\,\uqs$-covariantes. La sph\`ere quantique de  Podl\`es \cite{P}
est en effet une autre pr\'esentation de telles alg\`ebres, et
munies d'une involution. Nous n'avons pas besoin dans la suite
d'une quelconque involution sur l'alg\`ebre
$\,\displaystyle\,{\cal A}_{0,q}^c\,$.

$\bullet$ Dans la famille $\,\displaystyle\,{\cal
A}_{\hbar,q}^c\,$, nous traitons l'alg\`ebre
$\,\displaystyle\,{\cal A}_{0,q}^c\,$ comme la $q$-analogue de
l'alg\`ebre commutative $\,\displaystyle\,{\cal A}_{0,1}^c$.

\section{$\,sl(2)\,$ tress\'e}

\subsection{Crochet de Lie tress\'e de $sl(2)$ }

Le crochet de Lie $\,\displaystyle\,\lbc , \rbc\,$ sur $\,sl(2)$ :
\begin{eqnarray*}
\displaystyle\,\lbc ,\rbc\,:\,{\sf V}^{\ot 2}\,\to \,{\sf V},
\end{eqnarray*}

\noindent $\bullet\,\,$ est une application
$\,{\bbbk}$-lin\'eaire,

\noindent $\bullet\,\,$ et $\,sl(2)$-invariante. \vskip 0.5truecm

Nous allons d\'efinir de fa\c con analogue le crochet de Lie
tress\'e (not\'e $\,\displaystyle\,{\lbc ,\rbc}_q$) sur $\,sl(2)$.

Nous introduisons les $q$-analogues not\'ees $\,I_{\pm}^q\,$ des
sous-espaces sym\'e\-triques et antisym\'etriques de l'espace
$\,sl(2)^{\ot 2}\,$ de fa\c con similaire au cas classique en
posant : $$I_+^q\,=\,{\sf V}_0^q \op {\sf V}_2^q \quad
\mbox{et}\quad I_-^q\,=\,{\sf V}_1^q.$$ (Notons que les
alg\`ebres correspondantes $\,\displaystyle\,T({\sf
V})/\{I_{\pm}^q\}\,$ sont des d\'eforma\-tions plates de leurs
analogues classiques. En outre comme l'alg\`ebre
$\,\displaystyle\,{\cal A}_{0,q}^c\,$ est une alg\`ebre quotient
de l'alg\`ebre \lq\lq{$q$-sym\'etrique}" $\,\displaystyle\,T({\sf
V})/\{I_-^q\}\,$, $\,\displaystyle\,{\cal A}_{0,q}^c\,$ est
\'egalement une alg\`ebre $q$-sym\'etrique.)

Notons que  $\,\displaystyle\,I_-^q\,$ est engendr\'e par les
trois tenseurs
\begin{eqnarray*}
&\displaystyle\,q^2uv-vu,\quad
(q^3+q)(uw-wu)+(1-q^2)vv,\quad-q^2vw+wv.
\end{eqnarray*}

\begin{defin}\label{pup}
Le crochet de Lie tress\'e de $\,sl(2)\,$, est  l'op\'erateur
\begin{eqnarray*}
\displaystyle\,{\lbc ,\rbc}_q\, :\,{\sf V}^{\ot 2}\,\to \,{\sf
V}\,\,\,\, \mbox{v\'erifiant}\,\,
\end{eqnarray*}
\begin{enumerate}
\item $\,\displaystyle\,{\lbc ,\rbc}_q\,{I_+^q}\,=\,0\,$,
\item $\bullet\,\,$  $\,\displaystyle\,{\lbc ,\rbc}_q\,(q^2 uv-vu)\,=\,-\tau u$,

\noindent $\bullet\,\,$ $\,\displaystyle\,{\lbc ,\rbc}_q\,((q^3
+q)(uw-wu)+(1-q^2 )vv)\,=\,\tau v$,

\noindent $\bullet\,\,$ $\,\displaystyle\,{\lbc ,\rbc}_q\,(-q^2
vw+wv)\,=\,\tau w$.

$\,\displaystyle\,{\tau}\,$ est une constante non nulle.
\end{enumerate}
$\,\displaystyle\,{\sf V}\,$ muni du crochet
$\,\displaystyle\,{\lbc ,\rbc}_q\,$ est appel\'e {\it alg\`ebre
de Lie tress\'ee}  et not\'e $\,\displaystyle\,sl(2)_q$. Plus
pr\'ecis\'ement $$\,\,\,\displaystyle\,sl(2)_q :\,=\,({\sf
V}\,,\,{\lbc ,\rbc}_q).$$
\end{defin}
De fait ce crochet de Lie tress\'e d\'epend du facteur $\,\tau$.
Mais nous n\'egli\-geons cette d\'ependance en supposant que
$\,\tau\,$ est fix\'e.

\begin{prop}
\begin{enumerate}
\item $\,\displaystyle\,{\lbc ,\rbc}_q\,$ est un $\,\uqs$-morphisme (i.e.
une application $\,\uqs$-covariante).
\item La table de commutation de $\,\displaystyle\,{\lbc ,\rbc}_q\,$ est :
\begin{eqnarray*}
&\displaystyle\,{\lbc u,u\rbc}_q =0,\ \ {\lbc u,v\rbc}_q =-q^2
Mu,\
\ {\lbc u,w\rbc}_q =(q+q^{-1})^{-1}Mv,\\
&\displaystyle\,{\lbc v,u\rbc}_q =Mu,\ \ {\lbc v,v\rbc}_q =(1-q^2
)Mv,\
\ {\lbc v,w\rbc}_q =-q^2 Mw,\\
&\displaystyle\,{\lbc w,u\rbc}_q =-(q+q^{-1})^{-1} Mv,\
\ {\lbc w,v\rbc}_q =Mw,\ \ {\lbc w,w\rbc}_q = 0,\\
&\displaystyle\,M =(1+q^4 )^{-1}\tau .
\end{eqnarray*}
\end{enumerate}
\end{prop}

$Preuve :$ 1. C'est la propri\'et\'e 2 de la d\'efinition
\ref{pup}

2. C'est un calcul direct. \vskip 0.2truecm

(Si $\,q=1\,$ et $\,\displaystyle\,\tau =4\,$ (donc $\,M = 2\,$),
nous obtenons le crochet de Lie  sur $\,sl(2)$.)

\vskip 0.5truecm \noindent {\bf Remarque 3.1} C'est le fait que
l'espace $\,sl(2)\,$ apparaisse une seule fois dans la
d\'ecomposition en $\,sl(2)$-modules irr\'eductibles de
$\,sl(2)^{\ot 2}\,$ qui a  permit de d\'efinir de fa\c con unique
le crochet de Lie tress\'e $\,\displaystyle\,[,]_q$. Pour les
alg\`ebres de Lie
$\,\displaystyle\,{\sfg}=sl(n)\,\,(\,n\,>\,2\,)\,$, la
multiplicit\'e de l'espace $\,sl(n)\,$ dans $\,sl(n)^{\ot 2}\,$
est deux : l'une appartenant \`a la partie sym\'etrique de
$\,sl(n)^{\ot 2}\,$ et l'autre \`a sa partie antisym\'etrique. Il
n'est donc pas \'evident de d\'ecrire les  $q$-analogue des
alg\`ebres sym\'etriques et antisym\'etriques de $\,\sfg$.
Cependant, il existe un sous-espace
$\,\displaystyle\,I_-^q\subset {\sfg}_q^{\ot 2}\,$ o\`u
$\,{\sfg}_q\,$ est l'espace $\,sl(n)\,$ muni de la
$\,\displaystyle\,U_q(sl(n))$-action, telle que l'alg\`ebre
quadratique $\,\displaystyle\,T({\sfg}_q)/\{I_-^q\}\,$ soit une
d\'eformation plate de l'alg\`ebre sym\'etrique de $\,{\sfg}\,$
(voir \cite{D}). Une description explicite du sous-espace
$\,I_-^q\,$ peut \^etre donn\'ee par l'\'equation appel\'ee
\lq\lq{reflection equation}"
\begin{equation}
S L_1 S L_1\,=\,L_1 S L_1 S \label{ree}
\end{equation}
o\`u $\,S\,$ est une solution de l'\'equation de Yang-Baxter
quantique de type Hecke (voir \cite{G1}),
 $\,L_1 = L\ot id\,$ et $\,L\,$ est une matrice dont les
coefficients  matriciels sont les \'el\'ements $\,l_i^j, \quad
1\leq i\,,\,j\leq n$. L'alg\`ebre quadratique d\'efinie par
l'\'equation  (\ref{ree}) est habituellement appel\'ee
\lq\lq{{reflection equation algebra}" (REA). En consid\'erant
l'alg\`ebre $\,{\sfg}_q\,$ introduite dans  \cite{LS} qui est la
$q$-analogue de l'alg\`ebre $\,\sfg\,$, on peut d\'ecrire \`a
partir de la REA, l'alg\`ebre enveloppante de l'alg\`ebre
$\,{\sfg}_q\,$  qui est aussi une d\'eformation plate  \`a deux
param\`etres de l'alg\`ebre sym\'etrique de $\,\sfg$. (voir
\cite{AG})

Si $\,\sfg\,$ est une alg\`ebre de Lie simple diff\'erente de
$\,sl(n)\,$, dans $\,{\sfg}^{\ot 2}\,$ toute composante qui
appara\^\i t est sans multiplicit\'e. Donc on peut d\'efinir la
$q$-analogue du crochet de Lie, en imposant qu'il soit un
morphisme non trivial dans la cat\'egorie des
$\,U_q({\sfg})$-modules irr\'eductibles de dimension finie (il
est donc d\'efini de fa\c con unique \`a un facteur constant
pr\`es) puis, on introduit son alg\`ebre enveloppante, son
alg\`ebre sym\'etrique et antisym\'etrique comme dans le cas
classique, mais dans la cat\'egorie sous consid\'eration. (Ici le
fait que $\,{\sfg}^{\ot 2}\,$ soit sans multiplicit\'e joue le
r\^ole principale.)
 Cependant, ces alg\`ebres ne sont pas des d\'eformations plates de leurs
analogues classiques (voir \cite{G2}).

\subsection{Alg\`ebre enveloppante tress\'ee de $sl(2)$}

$\bullet\,$  Rappelons que dans le cas classique ($q=1$),
l'alg\`ebre enveloppante $\,\us\,$ est d\'efinie de la fa\c con
suivante
\begin{equation}
\begin{array}{ccc}
\displaystyle\,\us &=&T(sl(2))/\{ AB-BA-\lbc A,B\rbc \}\\
\displaystyle\,&=&T(sl(2))/\{ Im (id-{1\over 2}\lbc ,\rbc) I_-\},
\end{array}\label{aet}
\end{equation}
o\`u $\,A,\,B\,$ sont des \'el\'ements de $\,sl(2)\,$.

\noindent $\bullet\,$  Par analogie avec l'alg\`ebre
enveloppante  de $\,sl(2)\,$, nous d\'efinissons l'alg\`ebre
enveloppante tress\'ee de $\,sl(2)\,$ not\'ee
$\,\displaystyle\,\uovl\,$  comme suit :
$$\,\displaystyle\,\uovl\,=\,T(\ovl)/\{Im\,(id-\kappa\,{\lbc ,
\rbc}_q)I_-^q \}.$$ L'id\'eal
$\,\displaystyle\,\{Im\,(id-\kappa\,{\lbc ,\rbc}_q)I_-^q \}\,$
est celui engendr\'e par les \'el\'ements :
\begin{equation}
\begin{array}{cc}
&\displaystyle\,q^2 uv-vu-\kappa (q^2 {\lbc u,v\rbc}_q - {\lbc v,u\rbc}_q),\\
&\displaystyle\,(q^3 +q)(uw-wu)+(1-q^2 )vv-\kappa ((q^3 +q)({\lbc u,w\rbc}_q -\\
& -{\lbc w,u\rbc}_q )(1-q^2 ){\lbc v,v\rbc}_q),\\
&\displaystyle\,-q^2 vw+wv-\kappa (-q^2 {\lbc v,w\rbc}_q +{\lbc
w,v\rbc}_q ).
\end{array}\label{rat}
\end{equation}
Le choix de $\,\displaystyle\,\kappa\,$ sera  pr\'ecis\'e  par la
suite (voir la relation (\ref{vak})).

\vskip 0.5truecm \noindent {\bf Remarque 3.2} En fait on a
$${\cal A}_{\hbar,q}^c\,=\,\uovl/\{ {\cal C}_q -c \} \,\,
\mbox{avec}\,\,\hbar\,=\,\frac{\kappa \tau}{2}.$$

\begin{lemme}(\cite{DG})
Le Casimir tress\'e $\,\displaystyle\,{\cal C}_q\,$ est un
\'el\'ement centrale de l'alg\`ebre enveloppante tress\'ee
$\,\displaystyle\,\uovl\,$ c'est-\`a-dire : $$ X\,{\cal C}_q
\,=\,{\cal C}_q \qquad \forall\, \,X\,\in \,\,\uovl.$$
\end{lemme}

\vspace{0.5cm} Rappelons que nous voulons en fait d'une part,
introduire le module tangent sur l'hyperbolo\"\i de quantique et
d'autre part d\'ecrire les analogues tress\'ees de certaines
notions de la  g\'eom\'etrie diff\'erentielle classique sur ce
module. Pour ce faire, dans toute la suite  nous travaillons
principalement avec l'alg\`ebre $\,\displaystyle\,{\cal
A}_{0,q}^c\,$ et tous les modules consid\'er\'es sont des modules
sur cette  alg\`ebre.

\section{Espace tangent quantique}

Dans cette section, nous d\'efinissons l'espace tangent
(consid\'er\'e comme un $\,{\cal A}_{0,q}^c$-module) sur
l'hyperbolo\"\i de quantique not\'e $\,T({\rm H}_q)$. Puis nous
construisons des champs de vecteurs tress\'es sur l'hyperbolo\"\i
de quantique tels qu'ils nous permettent de munir $\,T({\rm
H}_q)\,$ d'une structure d'ancre quantique. \vskip 0.2truecm

Pour mieux pr\'esenter le formalisme de la construction de ce
module et de cette ancre dans le cas quantique, commen\c cons par
l'exemple du cas classique. \vskip 0.2truecm

Consid\'erons pour cela la sph\`ere $\,S^2\,$ de dimension deux.
Elle a pour \'equation $$x^2 + y^2 + z^2 = R^2$$ o\`u $\,R\,$ est
une constante strictement positive. Posons
$$\Fun(S^2)\,=\,{\bbbk}[x,\,y,\,z]/\{ x^2 + y^2 + z^2 - R^2\}.$$
Nous donnons ici, trois descriptions (globales) de l'espace
tangent sur la sph\`ere not\'e $\,T(S^2)$.

a) Comme un champ de vecteurs c'est : $\,\Vect(S^2)\,$ i.e.
l'espace des champs de vecteurs sur $\,S^2$ .

b) Comme un $\,\Fun(S^2)$-module : d'abord $\,\Vect(S^2)\,$ est
engendr\'e par les trois rotations infinit\'esimales
$\,X,\,Y,\,Z\,$ d\'efinies par (\ref{eqc}) et qui v\'erifient
dans l'alg\`ebre $\,\displaystyle\,\Fun(S^2)\,$ la relation
\begin{equation}
xX+yY+zZ\,=0.\label{eqo}
\end{equation}
(Dans l'expression (\ref{eqo}) $\,x,\,y,\,z\,$ d\'esignent en
fait des op\'erateurs de multiplication). Donc comme un
$\,\Fun(S^2)$-module, $\,T(S^2)\,$ peut \^etre r\'ealis\'e comme
le module quotient $\,M/N\,$ o\`u $$M\,=\,\{aX + bY + cZ\,,
\,a,\,b,\, c\,\in \,\Fun(S^2)\,\}\,,$$ $$N\,=\,\{f(xX+yY+zZ)\,,
f\,\in \,\Fun(S^2)\,\}\,.$$

c) Comme une vari\'et\'e alg\'ebrique affine : elle est plong\'ee
dans l'espace de dimension 6
$$(span\,(x,\,y,\,z,\,X,\,Y,\,Z\,))^*$$ et d\'efinie par
l'\'equation de la sph\`ere et la relation (\ref{eqo}). \vskip
0.2truecm

Notons que l'espace tangent $\,T(S^2)\,$ est un cas particulier
de fibr\'e vectoriel (sur la sph\`ere). Habituellement
$\,T(S^2)\,$ est d\'efini en termes de cartes locales. Nous
n'utilisons pas cette description locale ici, car dans le cas
quantique nous n'avons pas de localisation. \vskip 0.2truecm

En passant \`a l'hyperbolo\"\i de $\,{\rm H}\,$ (i.e. \`a
l'analogue non compacte de la sph\`ere),  nous avons de fa\c con
analogue trois descriptions (globales) de l'espace tangent sur
l'hyperbolo\"\i de not\'e $\,T({\rm H})$.

a') Comme un champ de vecteurs c'est juste $\,\Vect({\rm H})$.

b') Comme un $\,{\cal A}_{0,1}^c$-module, il est r\'ealis\'e
comme le module quotient $\,M/N\,$ o\`u $$M\,=\,\{aU + bV +
cW\,,\, a,\,b,\, c\,\in {\cal A}_{0,1}^c\,\,\}\,,$$
$$N\,=\,\{f(2uW + vV + 2wU)\,, f\,\in \,{\cal
A}_{0,1}^c\,=\,\Fun({\rm H})\}\,.$$ Notons qu'ici $\,U,\,V,\,W\,$
sont les rotations hyperboliques infinit\'esimales  associ\'ees
respectivement aux g\'en\'erateurs $\,u,\,v,\,w\,$ de l'alg\`ebre
de Lie $\,sl(2)$. De m\^eme l'analogue non compacte de
l'\'equation (\ref{eqo}) s'\'ecrit :
\begin{equation}
2u\,W + v\,V +2w\,U\,=\,0.\label{eqh}
\end{equation}

Nous pouvons mettre l'\'equation  (\ref{eqh}) sous la forme
symbolique
\begin{equation}
({\sf V}\ot {\sf V}')_0 = 0,\label{esh}
\end{equation}
o\`u $\,{\sf V}\,$ d\'esigne (encore) l'espace vectoriel
engendr\'e par $\,u,\,v,\,w\,$ et la marque $\,'\,$ d\'esigne
l'espace vectoriel engendr\'e par les rotations hyperboliques
infinit\'esimales. En outre  $\,{\sf V}\,$ et $\,{\sf V}'\,$ sont
des $\,\us$-modules. La composante $\,\displaystyle\,({\sf V}\ot
{\sf V}')_i\,$ d\'enote celle de spin $\,i\,$ dans la
d\'ecomposition en $\,\us$-modules irr\'eductibles de dimension
finie de $\,\displaystyle\,{\sf V}\ot {\sf V}'$.

Ici, nous regardons le module tangent $\,T({\rm H})\,$ comme un
$\,{\cal A}_{0,1}^c$-module gauche. Comme un $\,{\cal
A}_{0,1}^c$-module droit il est donn\'e par l'\'equation $\,({\sf
V}'\ot {\sf V})_0 =0$.  Il est bien connu que ces deux $\,{\cal
A}_{0,1}^c$-modules s'identifient naturellement. Nous discutons
dans la section 7, du probl\`eme d'identifica\-tion des modules
tangents gauche et droit  sur l'hyperbolo\"\i de quantique
(r\'ealis\'es tous deux comme des $\,{\cal A}_{0,q}^c$-modules).

c') Enfin comme une vari\'et\'e alg\'ebrique affine, elle est
plong\'ee dans l'espace de dimension 6
$$(span\,(u,\,v,\,w,\,U,\,V,\,W\,))^*$$ et d\'efinie par
l'\'equation de l'hyperbolo\"\i de et la relation  (\ref{eqh}).
\vskip 0.2truecm

Laquelle des trois descriptions pr\'ec\'edentes admet une
\lq\lq{bonne}" (i.e. plate) $q$-analogue ? \vskip 0.2truecm

Il est \'evident que si nous voulons d\'efinir sur
l'hyperbolo\"\i de quantique le $\,{\cal A}_{0,q}^c$-module
tangent gauche not\'e $\,T({\rm H}_q)_l\,$ comme une
d\'eformation plate de son analogue classique, nous devons
utiliser la m\^eme formule (\ref{esh}) mais dans la cat\'egorie
$\,\uqsm$. Cela  nous am\`ene \`a l'\'equation symbolique :
\begin{equation}
({\sf V}\ot {\sf V}'^q)_0 = 0.\label{sss}
\end{equation}
Dans l'identit\'e (\ref{sss}), $\,{\sf V}\,$ est muni de la
$\,\uqs$-action et $\,{\sf V}'^q\,$  est l'analogue tress\'e de
$\,{\sf V}'$. \vskip 0.2truecm

Interessons nous d'abord \`a la $q$-analogue des descriptions a')
et b').

\subsection{$T({\rm H}_q)\,$ comme $\,{\cal A}_{0,q}^c$-module tress\'e}

Soit $\,\displaystyle\,{\cal N}\,$ un \'el\'ement de la
cat\'egorie $\,\uqsm\,$ et qui soit en outre  un
$\,\displaystyle\,{\cal A}_{0,q}^c$-module. Soit l'application
$$\displaystyle\,\mu\,:\,{\cal A}_{0,q}^c \ot {\cal
N}\,\to\,{\cal N}$$ qui d\'esigne l'action de $\,{\cal
A}_{0,q}^c\,$ sur $\,{\cal N}$.
\begin{defin}\label{pipa}
$\,\displaystyle\,{\cal N}\,$ est appel\'e un
$\,\displaystyle\,{\cal A}_{0,q}^c$-{\it module tress\'e} si
$\,\displaystyle\,\mu\,$ est un morphisme dans la cat\'egorie
$\,\uqsm$ i.e. : $$z.\,\mu\,(a \ot n)\,=\,\mu\,(z_{(1)}.a \ot
z_{(2)}.n)\,$$ o\`u $\,\displaystyle\,z\,\in\,\uqs,\quad \Delta
(z)=z_{(1)} \ot z_{(2)},\,$ $\,a\in {\cal A}_{0,q}^c\,$ et
$\,\displaystyle\,n \in {\cal N}$.
\end{defin}

Notons $\,U^q,\,V^q,\,W^q\,$ les g\'en\'erateurs de $\,{\sf
V}'^q$. Le GQ $\,\uqs\,$ agit  sur ces g\'en\'erateurs comme il
op\'erait sur les \'el\'ements  de  l'espace $\,{\sf V}\,$ i.e.
\begin{eqnarray*}
&\displaystyle\,X.U^q\,=\,0,\,\,\,X.V^q=-(q+q^{-1})U^q,\,\,\,X.W^q=V^q,\\
&\displaystyle\,Y.U^q\,=\,-V^q,\,\,\,Y.V^q=(q+q^{-1})W^q,\,\,\,Y.W^q=0,\\
&\displaystyle\,H.U^q\,=\,2U^q,\,\,\,H.V^q=0,\,\,\,H.W^q=-2W^q.
\end{eqnarray*}
Alors l'identit\'e (\ref{sss}) devient en forme explicite :
\begin{equation}
(q^3 + q)uW^q + vV^q + (q + q^{-1})wU^q = 0 \label{casm}
\end{equation}

Ainsi le $\,{\cal A}_{0,q}^c$-module tangent gauche $\,T({\rm
H}_q)_l\,$ est
 r\'ealis\'e comme un $\,{\cal A}_{0,q}^c$-module facteur du
$\,{\cal A}_{0,q}^c$-module $$({\cal
A}_{0,q}^c)^3\,=\,M_l^q\,=\,\{aU^q + bV^q + cW^q,\,\, a,\,b,\,
c,\,\in {\cal A}_{0,q}^c\,\,\}\,$$ par le $\,{\cal A}_{0,q}^c$
sous-module $$N_l^q\,=\,\{f((q^3 + q)uW^q + vV^q + (q +
q^{-1})wU^q)\,, f\,\in \,{\cal A}_{0,q}^c\,\}.$$ Pr\'ecisons  que
les  $\,{\cal A}_{0,q}^c$-modules $$M_l^q,\quad {\rm N}_l^q,\quad
T({\rm H}_q)_l\,=\,M_l^q\,/\,{\rm N}_l^q\,$$ sont des modules
tress\'es au sens de la d\'efinition \ref{pipa}.

De fa\c con similaire le $\,{\cal A}_{0,q}^c$-module tangent
droit sur l'hyperbolo\"\i de quantique not\'e  $\,T({\rm
H}_q)_r\,$ est r\'ealis\'e comme le module quotient $\,
M_r^q/N_r^q\,$ o\`u $$M_r^q\,=\,\{{\ovu}^q\,a + {\ovv}^q\,b +
{\ovw}^q\,c\,,\,\, a,\,b,\, c,\,\in {\cal A}_{0,q}^c\,\,\}\,,$$
$$N_r^q\,=\,\{f((q^3 + q) {\ovu}^q\,w + {\ovv}^q\,v + (q + q^{-1})
{\ovw}^q\,u)\,,f\,\in \,{\cal A}_{0,q}^c\,\},$$ o\`u
$\,{\ovu}^q,\,{\ovv}^q,\,{\ovw}^q\,$ sont  les g\'en\'erateurs du
$\,{\cal A}_{0,q}^c$-module $\,T({\rm H}_q)_r$. (L'action du
groupe quantique $\,\uqs\,$ sur ces g\'en\'erateurs est la m\^eme
que sur  $\,{\sf V}'^q$.) De m\^eme les $\,{\cal
A}_{0,q}^c$-modules $$M_r^q,\quad {\rm N}_r^q,\quad T({\rm
H}_q)_r\,=\,M_r^q\,/\,{\rm N}_r^q\,$$ sont des modules tress\'es
au sens de la d\'efinition (\ref{pipa}).

\begin{prop}(\cite{A})
Le $\,{\cal A}_{0,q}^c$-module tangent $\,T({\rm H}_q)\,$ est une
d\'efor\-mation plate de son analogue classique.
\end{prop}

\subsection{$T({\rm H}_q)\,$ comme champs de vecteurs tress\'es }

Nous consid\'erons les g\'en\'erateurs $\,u,\,v,\,w\,$ de
l'alg\`ebre $\,{\cal A}_{0,q}^c\,$ comme des op\'erateurs de
multiplication (\`a gauche) dans cette m\^eme alg\`ebre, puis
nous d\'efinissons  les op\'erateurs $\,U^q,\,V^q,\,W^q\,$ qui
sont les $q$-analogues des rotations
 hyperboliques infinit\'esimales $\,U,\,V,\,W$. Finalement nous obtenons
la $q$-analogue des champs de vecteurs qui est engendr\'e par les
op\'erateurs $\,U^q,\,V^q,\,W^q\,$ et les \'el\'ements de
l'alg\`ebre  $\,{\cal A}_{0,q}^c\,$ sont trait\'es comme des
op\'erateurs. \vskip 0.2truecm

Soulignons une fois de plus que la tentation est de faire jouer
aux g\'en\'erateurs $\,X,\,Y,\,H\,$ du GQ $\,\uqs\,$ le r\^ole de
$q$-analogue des rotations hyperboliques infinit\'esimales
$\,U,\,V,\,W$. Mais alors dans ce cas, on a aucune relation de la
forme (\ref{sss}). Par cons\'equent ils ne peuvent  \^etre
consid\'erer comme les $q$-analogues  des champs de vecteurs
$\,U,\,V,\,W$. Il faut donc trouver un autre moyen pour
r\'ealiser les g\'en\'erateurs
$\,U^q,\,V^q,$\\
$\,W^q\,$ de $\,{\sf V}'^q\,$ comme des op\'erateurs sur
l'alg\`ebre $\,{\cal A}_{0,q}^c$.

Par analogie avec le cas classique, \`a partir du crochet de Lie
tress\'e de $\,sl(2)\,$  pr\'ec\'edemment d\'efini, on peut
d\'efinir $\,U^q,\,V^q,\,W^q\,$ comme  des op\'erateurs sur
l'espace  $\,\displaystyle\,{\sf V}\,$ (identifi\'e \`a l'espace
des fonctions  lin\'eaires sur $\,\displaystyle\,{\sf V}^*\,$).
En effet, associons au vecteur de base $\,u\,$ de l'espace
$\,\displaystyle\,{\sf V}\,$, l'op\'erateur $$
\begin{array}{cccl}
U^q :&{\sf V}&\to&{\sf V}\\
&z& \mapsto&\,\displaystyle\,U^q z\,=\,ad^q u(z)\,=\,[u,z]_q ,\ \
U^q 1=0,
\end{array}$$
d\'efini sur les vecteurs de  base de $\,\displaystyle\,{\sf V}$.
Les op\'erateurs $\,V^q,\,W^q\,$ associ\'es respectivement aux
vecteurs de base $\,v,\,w\,$ sont d\'efinis de fa\c con analogue
\`a l'op\'erateur $\,U^q$.

Imposons maintenant que les op\'erateurs  $\,U^q,\,V^q,\,W^q\,$
ainsi d\'efinis sur
 $\,\displaystyle\,{\sf V}\,$ v\'erifient les relations de d\'efinition de
l'alg\`ebre enveloppante tress\'ee $\,\uovl$. C'est-\`a-dire que
l'\'equation (\ref{rat}) soit encore satisfaite si nous rempla\c
cons les g\'en\'erateurs $\,u,\,v,\,w\,$ (de l'espace $\,{\sf
V}\,$) respectivement par leurs images par
$\,\,\displaystyle\,ad^q$. Autrement dit, pour l'op\'erateur
$\,\displaystyle\,ad^q\,$ nous avons pour  tout vecteur  de  base
$\,z\,$ de l'espace $\,\displaystyle\,{\sf V}$ :
\begin{equation}
\begin{array}{cc}
&q^2[u,[v,z]_q]_q-[v,[u,z]_q]_q=\kappa [q^2 [u,v]_q - [v,u]_q,z]_q \\
&\,\displaystyle\,(q^3 +q)([u,[w,z]_q]_q-[w,[u,z]_q]_q)+(1-q^2 )[v,[v,z]_q]_q=\\
&=\kappa (q^3 +q)[[u,w]_q-[w,u]_q,z]_q +(1-q^2 )[[v,v]_q,z]_q,\\
&\,\displaystyle\,-q^2[v,[w,z]_q]_q+[w,[v,z]_q]_q=\kappa  [-q^2
[v,w]_q +[w,v]_q,z]_q .
\end{array} \label{ijb}
\end{equation}
Les relations du (\ref{ijb}) fixent le choix de la constante
$\,\kappa\,$ de (\ref{rat}). Par exemple en se servant  de la
premi\`ere relation du (\ref{ijb}), nous obtenons
\begin{equation}
\,\displaystyle\,\kappa\,=1-(q^2+q^{-2})^{-1}.\label{vak}
\end{equation}
Ainsi pour ce choix de $\,\displaystyle\,\kappa\,$ fix\'e par
l'identit\'e (\ref{vak}), $\,\displaystyle\,ad^q\,$ est une
repr\'esenta\-tion de $\,\displaystyle\,\ovl\,$ (sur
$\,\displaystyle\,{\sf V}$). Nous regardons les relations du
(\ref{ijb}) comme {\it l'analogue tress\'e de l'identit\'e de
Jacobi}. \vskip 0.5truecm \noindent {\bf Remarque 4.2.1} Il a
\'et\'e montr\'e dans \cite{LS} qu'il
 existe \'egalement un analogue tress\'e de l'identit\'e de Jacobi pour les
alg\`ebres $\,sl(n),\,n>2$. Pour d'autres alg\`ebres de Lie
simples, il n'en existe (apparemment) pas. \vskip 0.2truecm

Puisque $\,ad^q\,$ est une repr\'esentation, les op\'erateurs
$\,\displaystyle\,U^q\,,\,V^q\,,\,W^q\,$ sont des op\'erateurs
adjoints (gauches) tress\'es. Ils sont donc d\'efinis par la
table suivante :
\begin{eqnarray*}
&\,\displaystyle\,U^q u=[u,u]_q=0,\ \ U^q v =-q^2 Mu,\ \
U^q w =(q+q^{-1})^{-1}Mv,\\
&\,\displaystyle\,V^q u
=Mu,\ \ V^q v =(1-q^2 )Mv,\ \ V^q w =-q^2 Mw,\\
&\,\displaystyle\,W^q u =-(q+q^{-1})^{-1} Mv,\ \ W^q v =Mw,\ \
W^q w = 0.
\end{eqnarray*}
\begin{lemme}
\begin{equation}
\displaystyle\,(q^3 +q)\,u\,W^q + v\,V^q +
(q+q^{-1})\,w\,U^q\,=\,0\,,\label{ean}
\end{equation}
sur les \'el\'ements de degr\'e un de l'alg\`ebre
$\,\displaystyle\,{\cal A}_{0,q}^c$.
\end{lemme}

$Preuve :$ Il suffit de la v\'erifier sur les \'el\'ements
$\,u,\,v,\,w\,$ de l'alg\`ebre $\,\displaystyle\,{\cal
A}_{0,q}^c$. Ce qui est imm\'ediat. Par exemple pour $\,u\,$ on a
: $$(q^3 +q)\,u\,W^q (u) + v\,V^q (u) + (q+q^{-1})\,w\,U^q
(u)\,=$$ $$=M(-q^2 uv + vu)\,=\,0\quad\mbox{dans
l'alg\`ebre}\quad {\cal A}_{0,q}^c.$$

\vskip 0.5truecm

Il reste maintenant \`a r\'esoudre le probl\`eme qui consiste \`a
\'etendre  les op\'erateurs $\,U^q,\,V^q,\,W^q\,$ (bien d\'efinis
pour le moment sur l'espace $\,{\sf V}\,$) sur tout \'el\'ement
de l'alg\`ebre $\,{\cal A}_{0,q}^c\,$ de telle mani\`ere que
l'\'equation (\ref{ean}) soit v\'erifi\'ee sur $\,{\cal
A}_{0,q}^c\,$ et les op\'erateurs ainsi prolong\'es satisfassent
aux relations de d\'efinition de l'alg\`ebre $\,\uovl$.

Dans le cas classique cette extension est effectu\'ee par la
r\`egle de Leibniz. Il existe aussi une forme de cette r\`egle
pour les solutions involutives de l'\'equation de Yang-Baxter
quantique (voir (\ref{san})). Mais dans notre cas,  o\`u
l'op\'erateur de Yang-Baxter quantique provenant du GQ $\,\uqs\,$
n'est pas involutif, ce n'est pas \'evident  d'\'etendre de fa\c
con naturelle ces op\'erateurs sur les \'el\'ements de degr\'e
sup\'erieur (\`a un) de l'alg\`ebre  $\,\displaystyle\,{\cal
A}_{0,q}^c$.

Le fait que l'alg\`ebre  $\,\displaystyle\,{\cal A}_{0,q}^c\,$
se  d\'ecompose en $\,\uqs$-modules irr\'educ\-tibles de spin
$\,k\,$, $\,{\sf V}_k^q (\subset {\sf V}^{\ot k})\,$ nous permet
de d\'efinir cette extension sur les composantes $\,{\sf V}_k^q$.

\begin{thm}
Il existe une application $$\beta\,:\, T({\rm H}_q) \ot {\cal
A}_{0,q}^c\,\to\,{\cal A}_{0,q}^c\,$$ telle que :

1) Le diagramme suivant soit commutatif
 $$
\begin{array}{ccc}
{\cal A}_{0,q}^c \ot T({\rm H}_q) \ot {\cal A}_{0,q}^c&\to &T({\rm H}_q)
\ot {\cal A}_{0,q}^c\\
\downarrow&&\downarrow\\
{\cal A}_{0,q}^c \ot {\cal A}_{0,q}^c&\to  &{\cal A}_{0,q}^c
\end{array}
$$ (o\`u dans les lignes horizontales de ce diagramme $\,{\cal
A}_{0,q}^c\,$ agit sur $\,{\cal A}_{0,q}^c\,$ par son produit
habituel et dans les lignes verticales on op\`ere par
l'application $\,\beta$).

2) L'application $\,\beta\,$ restreinte \`a l'espace $\,{\sf
V}'^q\,$ est une repr\'esentation de l'alg\`ebre $\,\uovl$. En
outre les g\'en\'erateurs de $\,\beta ({\sf V}'^q)\,$ v\'erifie
l'identit\'e (\ref{ean}).
\end{thm}

$Preuve :$

2) Elle se fait en  deux \'etapes :

{\bf Etape} 1. Nous allons r\'ealiser les op\'erateurs
$\,U^q,\,V^q,\,W^q\,$ comme une s\'erie de repr\'esentations de
l'alg\`ebre  $\,{\cal A}_{0,q}^c\,$ sur les composantes $\,{\sf
V}_k^q$.

Soit $$
\begin{array}{cccl}
\displaystyle\,P_k^q :&{\sf V}^{\ot k }&\to &{\sf V}_k^q\\
&z&\mapsto& P_k^q z
\end{array}$$
le projecteur tress\'e correspondant \`a son analogue classique
not\'e $\,\displaystyle\,P_k$. $\displaystyle\,P_k^q\,$  est un
morphisme dans la cat\'egorie $\,\uqsm$. Consid\'erons \'egalement
l'application
\begin{equation}
{\rho}_k^q (z) v\,=\,{\alpha}_k P_k^q ({\rho}^q (z) \ot
id_{k-1})v, \,\,z\,\in \,\ovl,\, v\in {\sf V}_k^q,\,{\alpha}_k
\in {\bbbk}\label{ero}
\end{equation}
o\`u $\,\displaystyle\,id_{k-1}\,$ est l'op\'erateur identit\'e
sur l'espace $\,\displaystyle\,{\sf V}^{\ot (k-1)}$. Notons que
dans la formule (\ref{ero}), $\,{\alpha}_k\,$ est une constante
quelconque. En supposant que l'application
$\,\displaystyle\,{\rho}_k^q\,$ soit une repr\'esentation de
l'alg\`ebre $\,{\cal A}_{0,q}^c\,$, cela impose le choix de la
constante $\,\displaystyle\,{\alpha}_k\,$ donn\'ee  par la
propositon suivante

\begin{prop}(\cite{G2})\label{alp}
$${\alpha}_k\,=\,(q^{-1}+ q^3)(1+ q^2 + ... + q^{2(k-1)})(q^{-1}+
q^{2k+1})^{-1}.$$
\end{prop}
\vskip 0.5truecm

Ainsi \`a l'\'el\'ement de base $\,u\,$ de l'espace $\,{\sf
V}\,$, on associe la famille d'op\'erateurs
$\,\displaystyle\,U_k^q\,$ d\'efinie sur le $\,\uqs$-module
irr\'eductible de spin $\,k\,$, $\,\displaystyle\,{\sf
V}_k^q\,\subset\,{\sf V}^{\ot k}$ : $$U_k^q\,:\,{\sf
V}_k^q\,\to\,{\sf V}_k^q\,$$ $$U_k^q\,(v)\,=\,{\alpha}_k\, P_k^q
(U^q \ot id_{k-1}) v, \quad v\,\in\,{\sf V}_k^q$$ o\`u
$\,\displaystyle\,{\alpha}_k\,$ est donn\'ee par la proposition
\ref{alp}. La famille d'op\'erateurs $\,U_k^q\,$ ($\,k\in
{\bbbn}\,$) d\'efinie l'extension (encore not\'ee $\,U^q$) de
l'op\'erateur $\,U^q\,$ (initialement d\'efinit sur $\sf V$)
 sur l'alg\`ebre $\,\displaystyle\,{\cal A}_{0,q}^c$
$$ U^q\,:\,{\cal A}_{0,q}^c\,\to\,{\cal A}_{0,q}^c.$$ On prolonge
de la m\^eme mani\`ere, les op\'erateurs
$\,\displaystyle\,\,V^q,\,\,W^q$. Nous montrons dans la deuxi\`eme
\'etape que (pour toute valeur de la constante  $\,{\alpha}_k\,$)
l'\'equation (\ref{ean}) est encore satisfaite par les
op\'erateurs prolong\'es $\,U^q,\,V^q,\,W^q$.

Pour le faire montrons d'abord  que l'\'el\'ement $\,(U^q \ot
id_{k-1})v_k^q\,$ (o\`u $\,v_k^q\,\in {\sf V}_k^q\,$), r\'eduit
en forme de base dans l'alg\`ebre $\,{\cal A}_{0,q}^c\,$ est
identique (\`a un facteur  pr\`es) \`a l'\'el\'ement $\,P_k^q(U^q
\ot id_{k-1})v_k^q$. Pour cela, il suffit  en fait de montrer que
dans la forme r\'eduite de l'\'el\'ement $\,(U^q \ot
id_{k-1})v_k^q\,$, il n'y a pas d'\'el\'ements appartenant aux
composantes $\,{\sf V}_i^q,\,i<k$. L'\'el\'ement $\,(U^q \ot
id_{k-1})v_k^q\,$ r\'esulte de l'application de l'op\'erateur
$\,[\,,\,]_q\,$ \`a l'\'el\'ement $\,u\ot v_k^q\,\in\,sl(2)_q \ot
{\sf V}_k^q$. Cette op\'eration commute avec l'action du GQ
$\,\uqs$.

Dans la d\'ecomposition en $\,\uqs$-modules du produit $\,sl(2)_q
\ot {\sf V}_k^q\,$ il y a trois composantes irr\'eductibles :
$\,{\sf V}_{k+1}^q,\,{\sf V}_k^q,\,{\sf V}_{k-1}^q$. En tenant
compte du fait que l'\'el\'ement $\,(U^q \ot id_{k-1})v_k^q\,\in
{\sf V}^{\ot k}\,$ et que dans la r\'eduction de base d'un
\'el\'ement de $\,{\sf V}^{\ot k}\,$ seulement les composantes
$\,{\sf V}_k^q,\,{\sf V}_{k-2}^q,\,{\sf V}_{k-4}^q\,$ peuvent
appara\^\i tre, nous en concluons que l'\'el\'ement $\,(U^q \ot
id_{k-1})v_k^q\,\in\,{\sf V}_k^q$. Dans le raisonnement
pr\'ec\'edent en rempla\c cant $\,U^q\,$ par $\,V^q\,$ ou par
$\,W^q\,$ on obtient la m\^eme conclusion. \vskip 0.5truecm

{\bf Etape} 2. Ayant finalement r\'ealis\'e $\,U^q,\,V^q,\,W^q\,$
comme des op\'era\-teurs sur l'alg\`ebre $\,\displaystyle\,{\cal
A}_{0,q}^c\,$, nous montrons \`a
 pr\'esent que l'\'equation (\ref{ean}) reste encore vraie sur toute l'alg\`ebre
$\,\displaystyle\,{\cal A}_{0,q}^c$. En effet, soit
$\,\displaystyle\,u_i\,\in\,\{u,\,v,\,w\}\,$ et
$\,\displaystyle\,g\,\in\,{\cal A}_{0,q}^c$. L'\'el\'ement
$\,\displaystyle\,g\,$ s'\'ecrit : $$g\,=\,\sum_k
\,P_k^q\,(g)\,=\,\sum_k g_k,\quad g_k\,=\,P_k^q\,(g) \quad
\mbox{avec} \quad g_k\,\in\,{\sf V}_k^q\,\subset\,{\sf V}^{\ot
k}.$$ Ecrivons \'egalement $\,g_k\,$ sous la forme $\,u_i
g_{k-1}$. Pour $\,k\,$ fix\'e, nous avons : $$\,[(q^3 +q)uW^q
+vV^q +(q+q^{-1})wU^q]\,(g_k)\,=\,(*)$$ $$={\alpha}_k[(q^3
+q)\,u\,P_k^q\,(W^q \ot id)+ v\,P_k^q\,(V^q\ot id) +
(q+q^{-1})w\,P_k^q\,(U^q\ot id)](g_k)$$ $\,(*)\,$ (\'ecrit dans
la base de  $\,{\cal A}_{0,q}^c\,$) donne (\`a un facteur
constant pr\`es) : $$(*)={\alpha}_k[(q^3 +q)\,P_k^q\,(uW^q \ot
id) + P_k^q\,(vV^q\ot id) + (q+q^{-1})P_k^q\,(wU^q\ot id)](g_k)$$
$$=\,{\alpha}_k P_k^q([(q^3 +q)uW^q u_i+vV^q u_i+(q+q^{-1})wU^q
u_i]\ot id_{k-1}(g_{k-1})\,=\,0.$$ La derni\`ere \'egalit\'e
\'etant due \`a l'\'equation (\ref{ean}). \vskip 0.2truecm

1) Elle est une cons\'equence imm\'ediate de la fa\c con dont les
op\'erateurs $\,U^q,\,V^q,\,W^q\,$ ont \'et\'e construits dans le
2). \vskip 0.5truecm

\begin{defin}\label{dcv}
Les op\'erateurs $\,\displaystyle\,U^q,\,V^q,\,W^q\,$ et toutes
leurs combinai\-sons lin\'eaires \`a coefficients dans l'alg\`ebre
$\,\displaystyle\,{\cal A}_{0,q}^c\,$ sont appel\'es les champs
de vecteurs tress\'es gauches. Ils sont de la forme $$a\,U^q +
b\,V^q + c\,W^q,\quad a,\, b,\, c\,\in\,{\cal A}_{0,q}^c.$$
\end{defin}
\vskip 0.5truecm {\bf Remarque 4.2.2} Pour  $\,q=1$, les
op\'erateurs $\,\displaystyle\,U^1,\,V^1,\,W^1\,$ co\"\i ncident
respectivement avec les champs de vecteurs (ou les rotations
hyperboliques infinit\'esimales)  $\,\displaystyle\,U,\,V,\,W\,$
sur l'hyperbolo\"\i de.

De la m\^eme mani\`ere, on d\'efinit
$\,\displaystyle\,{\ovu}^q,\,{\ovv}^q,\,{\ovw}^q\,$ les champs de
vecteurs tress\'es droits associ\'es respectivement aux
\'el\'ements de base $\,u,\,v,\,w\,$ de $\,{\sf V}$. Il est
facile de v\'erifier que dans l'alg\`ebre $\,\displaystyle\,{\cal
A}_{0,q}^c\,$, ces op\'erateurs satisfont \`a  une relation
analogue au (\ref{ean}) : $$ (q^3+q)\,{\ovu}^q\,w + {\ovv}^q\,v +
(q+q^{-1})\,{\ovw}^q u\,=\,0.$$ \vskip 0.5truecm

Par le th\'eor\`eme 4.1, $\,T({\rm H}_q)\,$ est le $\,{\cal
A}_{0,q}^c$-module tangent des champs de vecteurs tress\'es sur
l'hyperbolo\"\i de quantique. \vskip 0.5truecm \noindent {\bf
Remarque 4.2.3} Par la m\'ethode d\'evelopp\'ee dans \cite{LS}
qui consiste \`a d\'ecrire  l'alg\`ebre enveloppante tress\'ee
\`a partir de la REA, on aurait p\^u traiter facilement les
op\'erateurs adjoints tress\'es correspondants aux
g\'en\'erateurs $\,u,\,v,\,w\,$ comme  des op\'erateurs sur
l'alg\`ebre $\,{\cal A}_{0,q}^c$. Mais la diff\'erence
fondamentale avec notre m\'ethode est que les op\'erateurs
provenant de la m\'ethode sugg\'er\'ee dans \cite{LS}  ne
permettent  pas de contr\^o\-ler l'identit\'e (\ref{ean}) sur
l'alg\`ebre $\,{\cal A}_{0,q}^c$. \`A notre connaissance \`a part
notre fa\c con de d\'efinir les op\'erateurs
$\,\displaystyle\,U^q,\,V^q,\,W^q\,$ il n'en existe apparement
pas d'autre qui puisse contr\^oler (\ref{ean}).

\begin{defin}\label{ancq}
Le plongement $$\ovl\,\hookrightarrow\,T({\rm H}_q)\,$$ est
appel\'e une {\it ancre quantique}.
\end{defin}
Notons que le plongement de la d\'efinition \ref{ancq} est la
$q$-analogue de celui d\'efinit par (\ref{plon}), qui est
l'exemple le plus simple d'une ancre. Rappellons qu'une ancre est
constitu\'ee d'une vari\'et\'e $\,{\cal M}\,$ d'une alg\`ebre de
Lie $\,{\sf g}\,$ et d'un plongement de  $\,{\sf g}\,$ dans
l'espace des champs de vecteurs sur $\,{\cal M}$. C'est la raison
principale pour laquelle nous appellons le plongement de  la
d\'efinition (\ref{ancq}) ancre quantique et cela en d\'epit du
fait que le $\,{\cal A}_{0,q}^c$-module $\,T({\rm H}_q)\,$ n'est
muni d'aucun crochet de Lie tress\'e. Nous consid\'erons
\'egalement le couple ($\,T({\rm H}_q)\,,\,{\cal A}_{0,q}^c$)
comme une $q$-analogue d'une alg\`ebre de Lie-Rinehart \cite{R}
partielle (partielle est associ\'e au fait que $\,T({\rm H}_q)\,$
n'est muni d'aucune structure de Lie tress\'e). \vskip 0.5truecm

\noindent {\bf Remarque 4.2.4} Quant \`a la description c) qui
d\'ecrivait l'espace tangent sur l'hyperbolo\"\i de comme une
vari\'et\'e alg\'ebrique affine, pour la $q$-d\'eform\'e il est
n\'ec\'essaire de trouver la $q$-analogue de l'alg\`ebre
sym\'etrique de l'espace tangent $\,T({\rm H}_q)\,$ tel qu'il
soit une d\'eformation plate de son analogue classique. Le
probl\`eme fondamentale qui se pose pour l'exis\-tence de cette
alg\`ebre sym\'etrique d\'eform\'ee est de trouver une fa\c con
raisonnable de transposer les \'el\'ements de l'alg\`ebre
$\,{\cal A}_{0,q}^c\,$ et de l'espace $\,{\sf V}'^q$. Le seul bon
candidat susceptible de pouvoir r\'ealiser une telle
transposition est l'op\'erateur de Yang-Baxter (YB) quantique
provenant de la R-matrice universelle du GQ $\,\uqs$.
Malheureusement cette m\'ethode conduit \`a une d\'eformation non
plate de l'alg\`ebre sym\'etrique (classique). \vskip 0.5truecm
\noindent {\bf Remarque 4.2.5} Une fa\c con d'introduire dans le
cas classique les champs de vecteurs sur une vari\'et\'e
alg\'ebrique affine, consiste  \`a d\'efinir d'abord les champs
de vecteurs dans l'espace ambiant  comme toutes les combinaisons
lin\'eaires \`a coefficients-fonctions (ici les fonctions dans
l'espace ambiant) de d\'eriv\'ees partielles.  Ensuite, on
d\'efinit les champs de vecteurs sur la vari\'et\'e donn\'ee
comme de tels champs de vecteurs  qui respectent les \'equations
d\'efinissant la vari\'et\'e en question. L'autre fa\c con
consiste au passage aux cartes, i.e. au consid\'eration locale.

Malheureusement nous ne connaissons pas de $q$-analogues des
d\'eri\-v\'ees partielles. Par cons\'equent nous ne savons pas
d\'efinir les champs de vecteurs tress\'es sur l'espace
$\,{\ovl}^*\,$ tout entier, puisque les d\'eriv\'ees partielles
tress\'ees ne sont pas d\'efinies. C'est la raison principale
pour laquelle nous avons introduit les champs de vecteurs
tress\'es  sur l'hyperbolo\"\i de quantique \`a partir du crochet
de Lie tress\'e. \vskip 0.5truecm

Nous montrons maintenant que sur la sph\`ere (ou l'hyperbolo\"\i
de), la notion de champs de vecteurs d\'efinie \`a partir des
d\'eriv\'ees partielles et celle d\'efinie \`a partir des champs
de vecteurs adjoints sont \'equivalentes. \vskip 0.2truecm

En effet, il est clair que les champs de vecteurs $\,X,\,Y,\,Z\,$
sur la sph\`ere d\'efinis par (\ref{eqc}) tels que $$X\,(x^2 +
y^2 +z^2 -R^2)\,=\,Y\,(x^2 + y^2 +z^2 -R^2)\,=\, Z\,(x^2 + y^2
+z^2 -R^2)\,=\,0$$ s'expriment en fonction des d\'eriv\'ees
partielles $\,\partial_x,\,\partial_y,\,\partial_z$. Montrons la
r\'eciproque. Cela revient \`a montrer que  pour tout champ de
vecteurs $\,\displaystyle\,{\cal X}\,$ sur la sph\`ere de la forme
$${\cal X}\,=\,\alpha\,\partial_x + \beta\,\partial_y +
\gamma\,\partial_z,$$ o\`u
$\,\displaystyle\,\alpha\,,\,\beta\,,\,\gamma \in \Fun(S^2)\,$
sont tels que
\begin{equation}
{\cal X}(x^2 + y^2 +z^2 -R^2)=0\quad \mbox{i.e.}\quad \alpha\,x +
\beta\,y + \gamma\,z\,=\,0,\label{eqs}
\end{equation}
$\,\displaystyle\,{\cal X}\,$ est une combinaison lin\'eaire \`a
coefficients-fonctions (ici dans \\
$\,\Fun(S^2)\,$) des champs de vecteurs $\,X,\,Y,\,Z$. Notons que
la condition (\ref{eqs}) signifie que le champ de vecteurs
$\,\displaystyle\,{\cal X}\,$ est tangent \`a la sph\`ere.
\begin{prop}\label{peg}
Il existe $\,\displaystyle\,k,\,l,\,m\,\in\,\Fun(S^2)\,$ tels que
:
\begin{equation} {\cal X}\,=\,\alpha \partial_x+\beta \partial_y+\gamma
\partial_z = kX+lY+mZ\,.\label{eqq}
\end{equation}
\end{prop}

$Preuve :$ En  appliquant
 (\ref{eqq}) respectivement \`a $\,x,\,y,\,z\,$ les fonctions
$\,\displaystyle\,\alpha,\,\beta,\,\,\gamma\,$ se mettent sous la
forme
\begin{equation}
\alpha\,=\,-zl + ym,\,\,\beta\,=\,zk - xm,\,\,\gamma\,=\,-yk +
xl\,.\label{eqb}
\end{equation} Consid\'erons les fonctions $\,k,\,l,\,m\,$ d\'efinies par
\begin{equation}
\displaystyle\,k=\frac{1}{R^2}(\beta z-\gamma
y),\,\,l=\frac{1}{R^2}(-\alpha z+ \gamma x),\,\,
m=\frac{1}{R^2}(\alpha y-\beta x).\label{eug}
\end{equation}
Par l'\'equation (\ref{eqs}), il est facile de v\'erifier que
$\,(k,\,l,\,m)\,$ ainsi donn\'e   v\'erifie (\ref{eqq}). \vskip
0.5truecm

Par le changement de base : $$\,u\,=i(x+iy),\quad v\,={\sqrt
2}z,\quad w\,=-i(x-iy)\,,$$ on se ram\`ene au cas de
l'hyperbolo\"\i de. On peut remarquer qu'un champ de vecteurs
$${\cal X}\,=\,\alpha \partial_u + \beta \partial_v + \gamma
\partial_w $$ \lq\lq{respecte}" l'\'equation $$ 2uw +
\frac{v^2}{2}-c \,=0\,$$ (i.e. $\,\displaystyle\,2\alpha w +
\beta v + 2\gamma u =0$) si et seulement si
$\,\displaystyle\,{\cal X}\,$ peut \^etre pr\'esent\'e sous la
forme : $${\cal X}\,=\,k\,U +l\,V + m\,W$$ o\`u
$\,\displaystyle\,U,\,V,\,W\,$ sont les rotations hyperboliques
infinit\'esimales.

\section{Projectivit\'e du $\,{\cal A}_{0,q}^c$-module $\,T({\rm H}_q)$}

Il existe plusieurs d\'efinitions \'equivalentes d'un module
projectif (voir par exemple \cite{L}). En rempla\c cant
homomorphisme par homomorphisme dans la cat\'egorie $\,\uqsm\,$
et module par module tress\'e, nous adoptons ces d\'efinitions et
propri\'et\'es  pour notre cas tress\'e.

Pour montrer que le module tangent sur l'hyperbolo\"\i de
quantique est un module projectif, nous traitons juste le cas
classique, i.e. le cas de la sph\`ere puis nous en d\'eduisons
les r\'esultats dans notre cas tress\'e. \vskip 0.2truecm

Comme nous l'avons d\'eja fait remarquer, bien que le module
tangent
 $\displaystyle\,T(S^2)\,$  ne soit pas un $\,\Fun(S^2)$-module libre,
il est par contre un $\,\Fun(S^2)$-module projectif. En effet
posons ici $\,A = \Fun(S^2)$. Soit
($\,\displaystyle\,X,\,Y,\,Z\,$) la base du $\,A$-module
$\,A^3\,$ et $\,N\,$ le $\,A$-module de type fini engendr\'e par
l'\'el\'ement $$xX + yY + zZ\,:=\,(x,y,z)\,$$ $\,N\,$ est  un
sous-module  de $\,A^3$.  Notons  $\,\displaystyle\,\ovn\,$ le
sous-module de $\,A^3\,$ engendr\'e  par les \'el\'ements
$$yX-xY\,:=\,(y,-x,0),\quad zY-yZ\,:=\,(0,z,-y)\,,$$
$$\,xZ-zX\,:=\,(-z,0,x).$$

\begin{prop}\label{eig}
En sens de $\,A$-module on a : $$A^3\,=\,N \op \ovn$$
\end{prop}

$Preuve :$ Notons $\,{\rm Q}\,$ le projecteur de $\,A^3\,$ tel que
$\,\displaystyle\,\Im {\rm Q}=N\,$ et d\'efini par :
$\,\forall\,\,(f,\,g,\,h)\,\in \,A^3\,$
\begin{eqnarray*}
{\rm Q}(f,\,g,\,h)&=&R^{-2}(fx + gy + hz)(x,\,y,\,z)\\
&=&R^{-2}(fx + gy + hz)(xX + yY + zZ).
\end{eqnarray*}

Montrer que l'intersection des deux $A$-modules $N$ et $\ovn$ est
r\'eduite \`a z\'ero, revient \`a montrer que
$\,\displaystyle\,\Ker\,{\rm Q}=\ovn$. Nous avons : $${\rm
Q}(y,-x,0)\,=\,{\rm Q}(0,z,-y)\,=\,{\rm Q}(-z,0,x)=0$$ par
cons\'equent $\,\displaystyle\,\ovn \subset \Ker {\rm Q}$. Il
reste \`a prouver que  $\,\displaystyle\,\Ker {\rm Q} \subset
\ovn$.

Soit  $\,\displaystyle\,{\cal X}\in \Ker {\rm Q}\,$, alors il
existe $\,\displaystyle\,(\alpha,\,\beta,\,\gamma) \in A^3\,$ tel
que $${\cal X}\,=\,\alpha\,X + \beta\,Y + \gamma\,Z \quad
\mbox{et}\quad {\rm Q}({\cal X})=0.$$ La deuxi\`eme condition
($\,\displaystyle\,{\rm Q}({\cal X})=0\,$) entra\^\i ne
$$\alpha\,x + \beta\,y + \gamma\,z\,=\,0.$$ Ainsi si
$\,\displaystyle\,{\cal X}\in \Ker {\rm Q}\,$, cela revient \`a
dire que $\,\displaystyle\,{\cal X}\,$ est de la forme : $${\cal
X}\,=\,\alpha\,X + \beta\,Y + \gamma\,Z \quad \mbox{avec}\quad
\alpha\,x + \beta\,y + \gamma\,z\,=\,0.$$ Donc montrer qu'un tel
champ $\,\displaystyle\,{\cal X}\,$ appartient \`a
$\,\displaystyle\,\ovn\,$ revient  \`a montrer la proposition
\ref{peg} (il suffit de remplacer dans cette proposition, les
d\'eriv\'ees partielles
$\,\partial_x,\,\partial_y,\,\partial_z\,$ respectivement  par
les champs $\,X,\,Y,\,Z\,$). Par cons\'equent on a bien
$\,\displaystyle\,\Ker {\rm Q} \subset \ovn$.

En outre comme $\,(f,\,g,\,h)\,=\,f\,(1,\,0,\,0) + g\,(0,\,1,\,0)
+ h\,(0,\,0,\,1),\,$ il est n\'ec\'essaire et suffisant de
conna\^\i tre cette d\'ecomposition pour les \'el\'ements
$\,X:=(1,\,0,\,0);\,Y:=(0,\,1,\,0);\,Z:=(0,\,0,\,1)\,$ en vue de
pouvoir d\'ecomposer tout \'el\'ement du $\,A$-module $\,A^3\,$
comme somme d'un \'el\'ement du $\,A$-module $\,N\,$ et du
$\,A$-module $\,\ovn$. Par un calcul directe nous avons :
\begin{eqnarray*}
{\rm Q}(X)&=&R^{-2}x(xX + yY + zZ)\\
{\rm Q}(Y)&=&R^{-2}y(xX + yY + zZ)\\
{\rm Q}(Z)&=&R^{-2}z(xX + yY + zZ).
\end{eqnarray*}
Consid\'erons en outre le projecteur not\'e  ${\rm P}\,$, de
$\,A^3\,$ sur le $\,A$-module   $\,\ovn\,$ : $${\rm P}\,=\,id -
{\rm Q}\,,$$ nous obtenons \'egalement par un calcul directe :
\begin{eqnarray*}
{\rm P}(X)&=&-R^{-2}(z(xZ - zX) - y (yX - xY))\\
{\rm P}(Y)&=&-R^{-2}(x(yX - xY) - z (zY - yZ))\\
{\rm P}(Z)&=&-R^{-2}(y(zY - yZ) - x (xZ - zX)).
\end{eqnarray*}
D'o\`u la preuve de la proposition. \vskip 0.2truecm

Nous avons alors les identifications suivantes :
$$A^3\,/\,\{N\}\,\approx\,\ovn,\quad A^3\,/\,\{ \ovn
\}\,\approx\,N.$$ Donc le $\,A$-module $\,\displaystyle\,T(S^2) =
A^3\,/\,\{N\}\,$ est r\'ealis\'e comme un sous-module
($\,\displaystyle\,\ovn\,$) de $\,A^3$
 ayant un module suppl\'ementaire ($N$).

Faisons enfin remarquer qu'en prenant $\,A=\Fun({\rm H}),\,$
$\,N\,$ le $A$-module engendr\'e par $\,2wU + vV + 2uW\,$ et
$\,\ovn\,$ celui engendr\'e par les g\'en\'erateurs de $\,({\sf
V}\ot {\sf V}')_1,\,$ la proposition \ref{eig} reste \'egalement
valable sur l'hyperbo\-lo\"\i de. \vskip 0.2truecm

Passons \`a pr\'esent au cas tress\'e. L'analogue tress\'e des
$\,\Fun({\rm H})$-modules $N$ et $\,\ovn\,$ sont  respectivement
le $\,{\cal A}_{0,q}^c$-module (disons gauche) $\,N_l^q\,$ et
$\,{\ovn}_l^q$. Rappelons  que $\,N_l^q\,$ est engendr\'e par
$$(q^3+q)uW^q + vV^q + (q+q^{-1})wU^q$$ et pr\'ecisons que
$\,{\ovn}_l^q\,$ est engendr\'e par les \'el\'ements : $$q^2 uV^q
- vU^q,\,\,(q^3+q)(uW^q - wU^q) + (1-q^2)vV^q,\,\, -q^2 vW^q +
wV^q.$$
\begin{prop}\label{pepe}
En sens de $\,{\cal A}_{0,q}^c$-module on a : $$\,({\cal
A}_{0,q}^c)^3=\,N_l^q \op {\ovn}_l^q$$
\end{prop}

$Preuve :$ Pour montrer que l'intersection de $\,N_l^q\,$ et
$\,{\ovn}_l^q\,$ est r\'eduite \`a z\'ero, notons  $\,{\rm
Q}_q\,$ l'analogue tress\'e du projecteur $\,{\rm Q}$. $\,{\rm
Q}_q\,$ est la projection de $\,({\cal A}_{0,q}^c)^3\,$ telle que
$\,\displaystyle\,\Im {\rm Q}_q = N_l^q\,$ et d\'efinie par
 $\,\forall\,\,g,\,h,\,k\,\in \,{\cal A}_{0,q}^c\,$
\begin{eqnarray*}
&{\rm Q}_q (g,\,h,\,k)=c^{-1}[gu + hv + kw] [(q^3+q)uW^q + vV^q
+(q+q^{-1})wU^q]
\end{eqnarray*}
Consid\'erons \lq\lq{l'analogue tress\'e de la proposition
\ref{peg}}" i.e. en rempla\c cant dans l'\'equation (\ref{eqq})
les d\'eriv\'ees partielles par les g\'en\'erateurs de $\,{\sf
V}'^q\,$ et les champs $\,X,\,Y,\,Z\,$ par les g\'en\'erateurs de
$\,({\sf V}\ot {\sf V}'^q)_1$. En outre l'analogue tress\'e de la
condition (\ref{eqs}) s'\'ecrit $$(\alpha\,U^q + \beta\,V^q +
\gamma\,W^q)\,({\cal C}_q - c)\,=\,0$$ avec
$\,\alpha,\,\beta,\,\gamma \in {\cal A}_{0,q}^c$. Alors les
coefficients $\,k,\,l,\,m\,\in {\cal A}_{0,q}^c\,$ de l'analogue
tress\'e de la proposition \ref{peg} sont donn\'es (au facteur
constant $\,c^{-1}\,$ pr\`es) par les g\'en\'erateurs de l'espace
$\,({\bf V}\ot {\sf V})_1\,$ o\`u l'espace $\,{\bf V}\,$ est
celui engendr\'e par les \'el\'ements $\,\alpha,\,\beta,\,\gamma$.

Le fait que $\,\Ker ({\rm Q}_q) = {\ovn}_l^q\,$ d\'ecoule (comme
dans le cas classique) de \lq\lq{l'analogue tress\'e de la
proposition \ref{peg}}".

Comme pour la sph\`ere (ou l'hyperbolo\"\i de classique) tout
\'el\'ement du $\,{\cal A}_{0,q}^c$-module $\,({\cal
A}_{0,q}^c)^3\,$ est la somme d'un \'el\'ement de  $\,N_l^q\,$ et
de $\,{\ovn}_l^q$. En effet nous avons :
\begin{eqnarray*}
{\rm Q}_q (U^q)&=&c^{-1} u ((q^3+q)uW^q + vV^q +(q+q^{-1})wU^q)\\
{\rm Q}_q (V^q)&=&c^{-1} v ((q^3+q)uW^q + vV^q +(q+q^{-1})wU^q)\\
{\rm Q}_q (W^q)&=&c^{-1} w ((q^3+q)uW^q + vV^q +(q+q^{-1})wU^q).
\end{eqnarray*}
En consid\'erant en outre l'analogue tress\'e not\'e $\,{\rm
P}_q\,$ du projecteur $\,{\rm P}\,$ qui est tel que $\,\Im({\rm
P}_q) = {\ovn}_l^q\,$ et \'evidemment d\'efinit par $${\rm
P}_q\,=\,id - {\rm Q}_q\,$$ nous obtenons alors par un calcul
directe :
\begin{eqnarray*}
{\rm P}_q (U^q)=&-q^{-2}c^{-1}\{q^2 u[(q^3+q)(uW^q - wU^q) + (1-q^2)vV^q]-\\
&-v(vU^q - q^2 uV^q)\}\\
{\rm P}_q (V^q)=&-q^{-2}c^{-1}\{-(q^3+q)u(-q^2 vW^q + wV^q) - (q^3+q)(q^2 uV^q -\\
&- vU^q)- (1-q^2)v[(q^3+q)(uW^q - wU^q) + (1-q^2)vV^q]\}\\
{\rm P}_q (W^q)=&-q^{-2}c^{-1}\{-q^2 v(q^2 vW^q - wV^q) - w[(q^3+q)(uW^q - wU^q) +\\
&+ (1-q^2)vV^q]\}
\end{eqnarray*}

Ainsi nous venons de prouver (comme dans le cas classique) que le
$\,{\cal A}_{0,q}^c$-module $\,T({\rm H}_q)_l\,$ est un $\,{\cal
A}_{0,q}^c$-module projectif. (De m\^eme le $\,{\cal
A}_{0,q}^c$-module $\,T({\rm H}_q)_r\,$ est projectif).

\section{M\'etrique et Connexion tress\'ees}

Nous d\'efinissons et montrons l'existence d'une
(pseudo)m\'etrique et d'une connexion (partiellement d\'efinie)
tress\'ees  sur le module tangent $\,T({\rm H}_q)$.

\subsection{(Pseudo)m\'etrique tress\'ee}
Rappelons que \lq\lq{pseudo}" signifie que son analogue classique
n'est pas d\'efinie  positive. Par la suite nous omettrons cette
pr\'ecision.
\begin{defin}
L'op\'erateur $$
\begin{array}{cccl}
<,>:&T({\rm H}_q)_l{\ot}_{\bbbk} T({\rm H}_q)_r&\to &
{\cal A}_{0,q}^c\\
&a\ot b&\mapsto&{\displaystyle <a\ot b>=<a,b> }
\end{array}$$
est une m\'etrique tress\'ee si : $\,\displaystyle\,\forall\,
P\in\,T({\rm H}_q)_l,\, Q\in\,T({\rm H}_q)_r,\,f\in\, {\cal
A}_{0,q}^c\,$
\begin{equation}
<fP\,,\,Q>\,=\,f<P\,,\,Q>,\quad <P\,,\,Qf>\,=\,<P\,,\,Q>f
\label{pext}
\end{equation}
et $\,\displaystyle\,<\,,\,>\,$ est $\,\uqs$-covariant,
c'est-\`a-dire $\,\displaystyle\,\forall z\in \uqs$ :
\begin{eqnarray*}
z.<a,b>&=&<\,,\,>\Delta (z).(a\ot b)=<\,,\,>(z_{(1)}\ot z_{(2)})(a\ot b)\\
&=&<z_{(1)}.a,z_{(2)}.b>\quad \mbox{avec}\quad (a\,,\,b)\,\in
\,T({\rm H}_q)_l \times T({\rm H}_q)_r.
\end{eqnarray*}
Si en outre on a :
\begin{equation}
<\,,\,>\,({\sf V}'^q\ot {\osf}'^q)_1\,=\,0,\label{pqs}
\end{equation}
la m\'etrique tress\'ee est dite $q$-sym\'etrique.
($\,{\osf}'^q\,$ co\"\i ncide avec $\,{\sf V}'^q\,$ comme un
espace vectoriel. Mais il engendre $\,T({\rm H}_q)\,$ comme un
$\,{\cal A}_{0,q}^c$-module droit.)
\end{defin}
Soulignons que pour le moment nous introduisons la m\'etrique
tress\'ee sur $\,T({\rm H}_q)_l{\ot}_{\bbbk} T({\rm H}_q)_r$.
Nous pourrons la prolong\'ee sur  $\,T({\rm
H}_q)_{\epsilon}{\ot}_{\bbbk} T({\rm H}_q)_{\epsilon}\,$
seulement apr\`es l'identification des modules tangents
$\,\displaystyle\,T({\rm H}_q)_l\,$ et $\,\displaystyle\,T({\rm
H}_q)_r\,$ (voir section \ref{idt}).
\begin{thm}
Il existe une unique (\`a un facteur pr\`es) m\'etrique
tress\'ee  et $q$-sym\'etrique sur le module tangent
$\,\displaystyle\,T({\rm H}_q)$. Cette m\'etrique tress\'ee
restreinte sur $\,\displaystyle\,{\sf V}'^q\ot {\osf}'^q\,$
fournit la table suivante :
$$<U^q\,,\,{\ovu}^q>=uu,\,\,<U^q\,,\,{\ovv}^q>=uv,\,\,
<V^q\,,\,{\ovu}^q>=vu,$$ $$<W^q\,,\,{\ovv}^q>=wv,\,\,
<V^q\,,\,{\ovv}^q>=(1-q^2 )vv-q^{-1}(1+q^2 )^2 uw\,,$$
$$<W^q\,,\,{\ovw}^q>=ww,\,\, <U^q\,,\,{\ovw}^q>= -q^{-1}(1+q^2
)^{-1} vv-q^2 uw\,,$$ $$<V^q\,,\,\ovw>=vw,\,\,
<W^q\,,\,{\ovu}^q>= -q(1+q^2 )^{-1} vv-q^{-2}wu.$$
\end{thm}

$Preuve :$ D\'ecrivons d'abord tous les couplages
$$<\,,\,>\,:\,{\sf V}'^q \ot {\osf}'^q\,\to \,{\cal A}_{0,q}^c\,$$
$\uqs$-covariants. Pour cela, nous d\'ecomposons $\,{\sf V}' \ot
{\osf}'\,$ en $\,\uqs$-modules irr\'eductibles de dimension
finie. La $\,\uqs$-covariance impose les conditions :
\begin{equation}
<\,,\,>({\sf V}'^q \ot {\osf}'^q)_2 = k{\sf V}_2^q,\quad
<\,,\,>({\sf V}'^q \ot {\osf}'^q)_0 = \gamma \label{puq}
\end{equation}
o\`u $\,k\,$ et $\,\gamma\,$ sont des constantes. (La
$q$-sym\'etrie exigera en outre la relation (\ref{pqs}).)

La seconde \'etape de la d\'emonstration consiste \`a
d\'eterminer les relations de d\'ependance entre les param\`etres
$k$ et $\,\gamma$. Nous la faisons par le biais des relations
caract\'erisant le module tangent tress\'e gauche et le module
tangent tress\'e droit. Cela revient \`a
 v\'erifier que la relation de d\'ependance d\'efinie par
\begin{equation}
<\,,\,>^{23}({\sf V}\ot {\sf V}'^q)_0 \ot {\osf}'^q = 0
\label{dkg}
\end{equation}
est compatible avec celle d\'efinie par
\begin{equation}
<\,,\,>^{12}{\sf V}'^q\ot ({\osf}'^q \ot {\sf V})_0 = 0.
\label{dka}
\end{equation}
\vskip 0.1truecm

On \'etend enfin ce couplage \`a $\,\displaystyle\,T({\rm
H}_q)_l{\ot}_{\bbbk} T({\rm H}_q)_r\,$ en utilisant la
propri\'et\'e (\ref{pext}). (Voir appendice {\bf A} pour la forme
explicite des \'equations symboliques).

Nous avons  d\'etermin\'es le couplage
$\,\displaystyle\,<\,,\,>\,$ avec la condition (\ref{dkg}). Mais
on  v\'erifie que la condition (\ref{dka}) est satisfaite pour ce
couplage ainsi d\'etermin\'e.

\subsection{Connexion tress\'ee}

Rappelons d'abord que dans le cas classique ($q=1$), une
connexion lin\'eaire sur l'espace des champs de vecteurs $E$
d'une vari\'et\'e alg\'ebrique r\'eguli\`ere (ou plus
g\'en\'eralement d'une vari\'et\'e lisse) $\,{\cal M}\,$ est
l'application not\'ee $\,\displaystyle\,\nabla\,$ $$
\begin{array}{cccl}
\nabla :&E \ot E&\to &E\\
&a\ot b &\mapsto &{\displaystyle \nabla_a b}
\end{array}$$
$\bbbk$-lin\'eaire et satisfaisant aux propri\'et\'es suivantes :
\begin{enumerate}
\item $\,\,\displaystyle\,\nabla_{fa}\,b=f\nabla_a\,b,\quad (a,b)\in E^2,\,\,
\quad f\in \,\mbox{Fun}({\cal M})$
\item  $\,\,\displaystyle\,\nabla_a\,{fb}=f\nabla_a\,b + (af)\,b\,$.
\end{enumerate}
Notons que la propri\'et\'e 2. est la r\`egle de d\'erivation de
Leibniz pour la  connexion $\,\displaystyle\,\nabla$.

Lorsque la  connexion $\,\nabla\,$ est  sans torsion (par exemple
connexion de Levi-Civita), on a en plus
$$[a\,,\,b]=\nabla_a\,b\,-\,\nabla_b\,a\,$$ (o\`u  $\,[,]\,$
d\'esigne le crochet de Lie des champs de vecteurs $\,a\,$ et
$\,b\,$), on en d\'eduit alors :
\begin{eqnarray*}
\nabla_a {fb}&=&\nabla_{fb} a - [fb ,a]\\
&=&\nabla_{fb} a - (f[b , a] - (af)b)\\
&=&\nabla_{fb} a -(\nabla_{fb} a -\nabla_{fa} b - (af)b)\\
&=&f\nabla_a b + (af)b.
\end{eqnarray*}
Ainsi si la connexion $\,\displaystyle\,\nabla\,$ est sans
torsion, la propri\'et\'e 2. d\'ecoule de la propri\'et\'e 1.
Ceci nous permet de nous passer de la r\`egle de d\'erivation de
Leibniz pour la connexion si cette derni\`ere est sans torsion.
\vskip 0.2truecm

Nous g\'en\'eralisons la notion de connexion lin\'eaire \`a notre
cas tress\'e. Mais nous construisons dans ce cas plut\^ot une
connexion partiellement d\'efinie sur le module tangent tress\'e.
(\lq\lq{Partiellement d\'efinie}" signifie qu'elle est d\'efinie
sur un sous-ensemble de $\,\displaystyle\,T({\rm H}_q) \ot T({\rm
H}_q)$).

\begin{defin}\label{pige}
L'op\'erateur $$
\begin{array}{cccl}
\nabla :&T({\rm H}_q)_l {\ot}_{\bbbk}{\sf V}'^q&
\to&T({\rm H}_q)_l\\
&a\ot b &\mapsto &{\displaystyle \nabla_a\,b}
\end{array}$$
est une \lq\lq{connexion}" tress\'ee  (sans torsion), s'il
v\'erifie les propri\'etŽ\'es :
\begin{enumerate}
\item $\,\nabla\,$ est un morphisme dans la cat\'egorie des $\uqs$-modules,
c'est-\`a-dire : $$z.\nabla_a b=\nabla_{z_{(1)}.a}\, z_{(2)}.b$$
$\,\displaystyle\, \forall\,a\in T({\rm H}_q)_l,\quad b\in {\sf
V}'^q,\quad z\in \uqs$.
\item $\,\,\displaystyle{\nabla_{fa}\,b=f\nabla_a\,b},\quad
\forall a\in T({\rm H}_q)_l,\quad b\in {\sf V}'^q,\quad f\in
{\cal A}_{0,q}^c$.
\item
\begin{equation}
a^{ij} \nabla_{X_i} Y_j = [X_i , Y_j],\quad
X_i\,,\,Y_j\,\in\,{\sf V}'^q\quad \mbox{avec}\quad a^{ij} X_i \ot
Y_j \in {\sf V}_1^q.\label{qto}
\end{equation}
\end{enumerate}
\end{defin}

Notons que (\ref{qto}) est la $\,q$-analogue de la notion de
connexion sans torsion. Si nous arrivons en outre \`a \'etendre
le $q$-crochet de Lie $\,[\,,\,]_q\,$ sur
$\,\displaystyle\,T({\rm H}_q)\,$ tout entier et \`a comprendre
l'alg\`ebre enveloppante de cette alg\`ebre de Lie
$q$-d\'eform\'ee (nous en avons besoin pour \'ecrire la partie
\`a gauche de la relation (\ref{qto})), nous pourrons prolonger
notre connexion (partiellement d\'efinie) sur
$\,\displaystyle\,T({\rm H}_q)_{\epsilon}{\ot}_{\bbbk} T({\rm
H}_q)_{\epsilon}\,$ en utilisant une analogue de  (\ref{qto}).
Malheureusement nous ne connaissons aucune fa\c con de le faire.

\begin{thm}
Il existe une connexion tress\'ee  $\,\nabla\,$ sur le module
tangent de l'hyperbolo\"\i de quantique (au sens de la
d\'efinition (\ref{pige}). Elle est donn\'ee sur
$\,\displaystyle\,{\sf V}'^q\,$ par
\begin{eqnarray*}
\nabla_{U^q} U^q=&\alpha (uvU^q-q^2 uuV^q),\quad \nabla_{W^q} W^q=
\alpha(wwV^q-q^4 vwW^q),\\
\nabla_{U^q} V^q=&\beta\{\{-\alpha[(q^3 +q)uw-q^2
vv]-2q^2\}U^q-\alpha(q^6 +q^2
-1)\\
&uvV^q \}+\beta\alpha(q^3 +q)uuW^q,\\
\nabla_{V^q} U^q=&\beta\{\{-\alpha q^2 [(q^3 +q)uw-q^2
vv]+2\}U^q-\alpha q^2
(q^6 +q^2 -1)\\
&uvV^q \}+\beta\alpha q^2 (q^3 +q)uuW^q,\\
\nabla_{V^q} W^q=&\beta\{-\alpha q^3 (1+q^2)wwU^q+\alpha(1+q^4
-q^6 )vwV^q \}+
\beta\{\alpha q^4\\
&[(q+q^{-1})uw-vv]-2q^2 \}W^q,\\
\nabla_{W^q} V^q=&\beta\{-\alpha q^5 (1+q^2 )wwU+\alpha q^2
(1+q^4 -q^6 )
vw^q \}+\beta\{\alpha q^6 \\
&[(q+q^{-1})uw-vv]+2\}W^q,\\
\nabla_{W^q} U^q=&\beta\{\alpha q^6 vwU^q+\{-q^4 (1-q^2 )\alpha
[-uw+[2]^{-1}v^2 ]-\frac{2}{1+q^2}\}\\
&V^q\}-q^6\beta \alpha uvW^q,\\
\nabla_{U^q} W^q=&\beta\{\alpha q^2 vwU^q+\{(q^2 -1 )\alpha
[-uw+[2]^{-1}v^2 ]+
\frac{2q}{1+q^2}\}V^q\\
&-q^2 \alpha uvW^q\},\\
\nabla_{V^q} V^q=&\beta\{-\alpha q^3 (1+q^2 )^2
vwU^q+\{q(1+q^2)(1-q^4 )
[-uw+[2]^{-1}\\
&v^2 ]+2(1-q^2 )\}V^q\}+
q^3 (1+q^2 )^2 \alpha uvW^q \},\\
&\alpha =-\frac{2}{(1-q^2 +q^4 )c},\quad \beta=(1+q^4 )^{-1}.
\end{eqnarray*}
\end{thm}

$Preuve :$ Comme dans le cas de la m\'etrique tress\'ee, nous
d\'ecrivons d'abord toutes les applications $${\sf
V}'^q{\ot}_{\bbbk}{\sf V}'^q\,\to\,T({\rm H}_q)_l$$
$\uqs$-covariantes. Cette $\,\uqs$-covariance impose les
conditions suivantes :
\begin{equation}
\nabla ({\sf V}'^q\ot {\sf V}'^q)_2 = \alpha ({\sf V}\ot {\sf
V}'^q)_2,\quad \nabla ({\sf V}'^q\ot {\sf V}'^q)_0 = 0,\,\alpha
\in {\bbbk} \label{cci}
\end{equation}
compl\'et\'ees par les relations qui d\'ecoulent de (\ref{qto}).

(Pr\'ecisons que $\,\displaystyle\,\nabla ({\sf V}'^q\ot {\sf
V}'^q)_2\,$ d\'esigne l'ensemble $\,\displaystyle\,\nabla_a b\,$
o\`u $\,a\ot b\,$ parcourent $\,({\sf V}'^q\ot {\sf V}'^q)_2$.)

Concr\`etement, $\,U^q \ot U^q\,$ \'etant de poids deux et
$\,\displaystyle\,X.\nabla_{U^q} U^q =0\,$, consid\'erons
\begin{equation}
\nabla_{U^q} U^q=\alpha (uvU^q-q^2 uuV^q) \label{nab}
\end{equation}
o\`u $\,uvU^q-q^2 uuV^q\,$ est de poids 2. Notons que le choix
 $\,\displaystyle\,\nabla_{U^q} U^q = uU^q\,$ n'est pas compatible avec la
relation $$\nabla_{K}= (q^3+q)\,u\nabla_{W^q}+ v\nabla_{V^q}+
(q+q^{-1})\, w\nabla_{U^q}=0 \,\,\mbox{dans}\,\,T({\rm H}_q)_l.$$
Notons $$J_i =\frac{1}{[i]_q}Y^{i} J_{i-1},\quad i\in
\{1,\,2,\,3,\,4\}, \quad [i]_q =\frac{q^i
-q^{-i}}{q-q^{-1}},\quad \mbox{o\`u}$$ $$J_0 = \nabla_{U^q} U^q =
\alpha (uvU^q-q^2 uuV^q).$$ Par application de l'op\'erateur
$\,\displaystyle\,Y^{i}\,$ \`a la relation (\ref{nab}) nous
obtenons les quatre relations suivantes : $$-\nabla_{U^q} V^q-q^2
\nabla_{V^q} U^q\,=J_1\,,$$ $$ -\nabla_{U^q} W^q+q\nabla_{V^q}
V^q-q^4\nabla_{W^q} U^q\,=J_2\,,$$ $$\nabla_{V^q}
W^q+q^2\nabla_{W^q} V^q\,=J_3\,,$$ $$ \nabla_{W^q} W^q\,=J_4.$$
(Voir appendice {\bf B} pour les formes explicites). Le choix de
$\,\alpha\,$ est pr\'ecis\'e par la deuxi\`eme \'equation du
(\ref{cci}). (voir appendice {\bf B} ) \vskip 0.5truecm

{\bf Remarque 6.2} \label{pipi} Bien que nous ne savons pas
d\'efinir la $\,q$-analogue de la notion de courbure (dans le cas
classique elle est introduite localement)  dans le cadre de notre
approche globale, on peut d\'eviner (\`a un facteur pr\`es) la
forme de courbure correspondante en supposant que cette forme
soit une d\'eformation plate de son analogue classique. Sur la
sph\`ere, cette forme est donn\'ee (\`a un facteur pr\`es) par
$$x\,(dy\,dz - dz\,dy) + y\,(dz\,dx - dx\,dz) + z\,(dx\,dy -
dy\,dx).$$ Ainsi dans le cas tress\'e, pour obtenir la forme de
courbure, il suffit de remplacer dans l'expression $\,(q^3 + q)uw
+ vv + (q + q^{-1})wu\,$ les facteurs de droites respectivement
par la $q$-analogue des formes $\,dv\,dw - dw\,dv,\,dw\,du -
du\,dw,\,du\,dv - dv\,du$.

\section{Identification de $\,T({\rm H}_q)_l\,$ et $\,T({\rm H}_q)_r$ }

Nous nous int\'eressons \`a pr\'esent au probl\`eme qui consiste
\`a identifier les modules tangents gauche $\,T({\rm H}_q)_l\,$
et droit $\,T({\rm H}_q)_r$. Cette identification est
n\'ecessaire, car elle permet d'\'etendre la m\'etrique
tress\'ee  sur $\,T({\rm H}_q)_l\,$ (ou $\,T({\rm H}_q)_r\,$).

Le fait que l'alg\`ebre $\,{\cal A}_{0,q}^c\,$ ne soit munie
d'aucune structure d'involu\-tion et que l'op\'erateur de tresse
provenant de la $R$-matrice universelle du GQ $\,\uqs$ ne soit
pas involutif d'autre part, nous ont amen\'ees \`a sugg\'erer
d'identifier les $\,{\cal A}_{0,q}^c$-modules $\,T({\rm
H}_q)_l\,$ et $\,T({\rm H}_q)_r\,$ de la fa\c con suivante
\noindent --  nous construisons une base de chacun des modules
tangents, \noindent -- puis nous construisons une application
(dans la cat\'egorie $\,\uqsm$) entre $\,T({\rm H}_q)_l\,$  et
$\,T({\rm H}_q)_r\,$ qui
 co\"\i ncide (pour $\,q=1\,$) avec l'applica\-tion d\'efinie
par la volte.

\subsection{Base du $\,{\cal A}_{0,q}^c$-module  $T({\rm H}_q)$}\label{idt}

Notons d'abord que la m\'ethode de construction de cette base
\'etant aussi bien valable dans le cas classique que dans notre
cas tress\'e, nous consid\'erons dans ce qui suit indiff\'erement
ces deux cas.

En consid\'erant les modules $${\sf V}\ot {\sf V}',\quad {\sf
V}_k \ot {\sf V}',\quad k=2,\,3,\,...$$ o\`u les $\,{\sf V}_k\,$
(composantes  de base de l'alg\`ebre $\,{\cal A}_{0,1}^c\,$) sont
des $\,\us$-modules, nous d\'eterminons les composantes qui
\lq\lq{survivent}" dans le module tangent. Il est \'evident que
dans le produit $\,{\sf V}\ot {\sf V}'\,$ seulement deux
composantes survivent. \`A savoir : $$({\sf V}\ot {\sf
V}')_1,\quad ({\sf V}\ot {\sf V}')_2$$ puisque par construction
la composante $\,({\sf V}\ot {\sf V}')_0\,$ est nulle dans le
module tangent.

De fa\c con analogue, dans le produit $\,{\sf V}_2 \ot {\sf
V}'\,$, les composantes $\,({\sf V}_2 \ot {\sf V}')_3,\,\, ({\sf
V}_2 \ot {\sf V}')_2\,$ \lq\lq{survivent}" et il n'est pas
difficile de montrer que la composante $\,({\sf V}_2 \ot {\sf
V}')_1\,$ est nulle modulo les termes de $\,{\bbbk}\ot {\sf V}' =
{\sf V}'$. En effet, par construction, les \'el\'ements de
$\,\mu^{12}({\sf V}\ot ({\sf V}\ot {\sf V}')_0)\,$ sont nuls dans
le module tangent. En r\'eduisant en outre tout \'el\'ement du
produit $\,{\sf V}\ot {\sf V}\,$ \`a sa forme canonique, il est
la  somme d'un \'el\'ement de $\,{\sf V}_2 \subset {\sf V}^{\ot
2}\,$ et d'un \'el\'ement de $\,{\bbbk}$. D'o\`u la preuve de ce
dernier fait.

Ainsi de fa\c con g\'en\'erale, dans le produit $\,{\sf V}_k \ot
{\sf V}'\,$ qui se d\'ecompose de la mani\`ere suivante $${\sf
V}_k \ot {\sf V}'=({\sf V}_k \ot {\sf V}')_{k+1}\op ({\sf V}_k
\ot {\sf V}')_k \op ({\sf V}_k \ot {\sf V}')_{k-1}\,,$$ la
composante $\,({\sf V}_k \ot {\sf V}')_{k-1}\,$ est nulle modulo
les termes appartenant \`a $\,{\sf V}_j \ot {\sf V}',\quad j<k-1$.

Nous avons donc montr\'e le fait suivant :

\begin{prop}\label{ome}
Une base dans le module tangent  gauche $\,T({\rm H}_q)_l\,$ est
form\'ee par les $\,\uqs$-modules $$1.\,{\sf V}'^q,\quad
2.\,({\sf V}\ot {\sf V}'^q)_{1,2},\quad 3.\,({\sf V}_2^q \ot {\sf
V}'^q)_{2,3},\quad 4.\,({\sf V}_3^q \ot {\sf V}'^q)_{3,4} \quad
\mbox{etc}.$$
\end{prop}

\vskip 0.1truecm

On construit de fa\c con similaire une base du module tangent
droit $\,T({\rm H}_q)_r$. \vskip 0.2truecm

Identifions \`a pr\'esent les modules tangents  $\,T({\rm
H}_q)_l\,$ et $\,T({\rm H}_q)_r\,$ en d\'efinissant :
$$\alpha\,:\,T({\rm H}_q)_l\,\to \,T({\rm H}_q)_r$$ telle que
l'application : $$\alpha\,:\,({\sf V}\ot {\sf
V}'^q)_1\,\to\,({\sf V}'^q\ot {\sf V})_1, \quad ({\sf V}\ot {\sf
V}'^q)_2\,\to \,({\sf V}'^q\ot {\sf V})_2\,,$$ $$({\sf V}_2^q\ot
{\sf V}'^q)_2\,\to \,({\sf V}'^q\ot {\sf V}_2^q)_2,...$$ soit un
$\,\uqs$-morphisme unique (\`a un facteur constant pr\`es sur
chaque composante). Nous sugg\'erons maintenant une fa\c con
\lq\lq{canonique}" d'\'eliminer ce degr\'e de libert\'e sur les
composantes, de la mani\`ere suivante. On identifie les
\'el\'ements de l'espace : $$({\sf V}\ot {\sf V}'^q)_2\quad
\mbox{et}\quad ({\sf V}'^q\ot {\sf V})_2,\quad ({\sf V}_k^q \ot
{\sf V}'^q)_{k+1}\quad \mbox{et}\quad ({\sf V}'^q\ot {\sf
V}_k^q)_{k+1},\,\, k=2\,,3,...$$ qui co\"\i ncident quand on
remplace $\,{\sf V}'^q\,$ par $\,{\sf V}$. Quant aux composantes
$$({\sf V}\ot {\sf V}'^q)_1\quad \mbox{et}\quad ({\sf V}'^q\ot
{\sf V})_1,\quad ({\sf V}_k^q \ot {\sf V}'^q)_k\quad
\mbox{et}\quad ({\sf V}'^q\ot {\sf V}_k^q)_k\,,$$ leurs
\'el\'ements sont identifi\'es si la m\^eme op\'eration am\`ene
aux images oppos\'ees. On peut facilement voir  que dans le cas
classique, cette identification et celle d\'efinie par la volte
co\"\i ncident : c'est la motivation de notre m\'ethode
d'identification canonique. \vskip 0.5truecm

{\bf Remarque 7.1} Lorsque l'alg\`ebre $\,{\cal A}_{0,q}^c\,$ est
munie d'un op\'erateur de tresse involutif ($\,S^2 = id\,$), une
telle identification est faite de fa\c con similaire au cas
classique en rempla\c cant la volte par $\,S$. Mais pour
l'op\'erateur de tresse non involutif provenant du GQ $\,\uqs\,$
ce n'est plus raisonnable. Consid\'erons $\,M_1\,$ et $\,M_2\,$
deux  $\,{\cal A}_{0,q}^c$-modules (gauche par exemple). Le
probl\`eme r\'eside dans le fait que dans le produit $$m_1 \ot
f\,m_2,\quad m_1 \in M_1,\, m_2 \in M_2,\,f \in {\cal
A}_{0,q}^c\,$$ il n'existe aucune fa\c con raisonnable de
transposer le facteur $\,f\,$ pour le mettre \`a gauche de telle
mani\`ere que le produit tensoriel $\,{\ot}_{{\cal A}_{0,q}^c}\,$
soit associatif et le module $\,,M_1{\ot}_{{\cal A}_{0,q}^c}
M_2\,$ soit une d\'eformation plate de son analogue classique en
supposant bien s\^ur que  $\,M_1\,$ et $\,M_2\,$ soient des
d\'eformations plates de leurs analogues classiques respectifs.

\vskip 0.5truecm

En r\'esum\'e, pr\'ecisons une fois de plus que notre m\'ethode
qui consiste \`a introduire la m\'etrique en deux \'etapes :

$\bullet \quad$ en d\'efinissant d'abord un couplage sur
$\,\displaystyle\,T({\rm H}_q)_l {\ot}_{\bbbk} T({\rm H}_q)_r\,$,

$\bullet \quad$ et ensuite \`a identifier
$\,\displaystyle\,T({\rm H}_q)_l\,$
et $\,\displaystyle\,T({\rm H}_q)_r\,$,\\
nous permet de contr\^oler le fait que notre construction ne soit
pas contradictoire (i.e. la platitude a lieu).

Il existe dans la litt\'erature plusieurs constructions de
fibr\'e  quantique sur la sph\`ere (\cite{BM}, \cite{S2}) et de
m\'etrique quantique. Mais pour la plupart d'entre elles, la
notion de platitude de ces constructions est carr\'ement
ignor\'ee. C'est en cela que nos constructions sont diff\'erentes.

\newpage

\section*{Appendice}

\subsection*{{\bf A} $\qquad $ Hyperbolo\"\i de quantique et
m\'etrique tress\'ee}

Nous avons explicitement :
\begin{eqnarray*}
&\displaystyle\,{\sf V}_0^q = Span\,((q^3+q)uw+vv+(q+q^{-1})wu ),\\
&\displaystyle\,{\sf V}_1^q =
Span\,(q^2uv-vu,\,\,(q^3+q)(uw-wu)+(1-q^2)vv,\,\,
-q^2vw+wv ),\\
&\displaystyle\,{\sf V}_2^q = Span\,
(uu,\,\,uv+q^2vu,\,\,uw-qvv+q^4wu,\,\,vw+q^2wv,\,\,ww).
\end{eqnarray*}
(Dans ces expressions nous avons omis le symbole du produit
tensoriel, ce que nous ferons lorsque cela n'apporte aucune
confusion.)

$\,\displaystyle\,{\cal C}_q\,$ est le g\'en\'erateur de l'espace
de dimension un  ($\,\displaystyle\,{\sf V}_0^q\,$) : $${\cal
C}_q\,=\,(q^3+q)uw+vv+(q+q^{-1})wu.$$

En posant : $$x_0\,=\,u,\quad x_1\,=\,q^2uv-vu,\quad x_2\,=\,uu$$
nous avons :
\begin{eqnarray*}
&\displaystyle\,q^2uv-vu=-2u\hbar,\,\,\,(q^3+q)(uw-wu)+(1-q^2)vv=2v\hbar,\,\,\\
&\displaystyle\,-q^2vw+wv=2w\hbar,\,\,\,{\cal
C}_q=(q^3+q)uw+vv+(q+q^{-1})wu=c,
\end{eqnarray*}

Explicitons \`a pr\'esent les \'equations symboliques de la
Section 6 :

$\bullet \quad <\,,\,>({\sf V}'^q \ot {\osf}'^q)_2 = k{\sf
V}_2^q$.

Sur la composante de spin 2 de $\,\displaystyle\,{\sf V}'^q\ot
{\osf}'^q\,,$ l'\'el\'ement $\,\displaystyle\,U^q\ot {\ovu}^q\,$
est de poids 2 et nous avons n\'ecessairement par la
$\,\uqs$-covariance : $$\,<U^q\,,\,{\ovu}^q>\,=\,kuu$$ Par
application de l'op\'erateur $\,Y^n\,$ (avec $\, n\in \{
1,\,2,\,3,\,4\}\,$) \`a $\,\displaystyle\,<U^q\,,\,{\ovu}^q>\,$
nous obtenons les quatre relations suivantes :
$$<U^q\,,\,{\ovv}^q>+ q^2<V^q\,,\,{\ovu}^q>= k(q^2vu+uv)\,,$$
$$-<U^q\,,\,{\ovw}^q>+ q<V^q\,,\,{\ovv}^q>-
q^4<W^q\,,\,{\ovu}^q>= k(qvv - q^4wu - uw)\,,$$
$$<V^q\,,\,{\ovw}^q> + q^2<W^q\,,\,{\ovv}^q>= k(vw + q^2wv)\,,$$
$$<W^q\,,\,{\ovw}^q> = kww.$$ $\bullet \quad <\,,\,>({\sf V}'^q
\ot {\osf}'^q)_0 = \gamma$.

Sur la composante de spin z\'ero de $\,\displaystyle\,{\sf
V}'^q\ot {\osf}'^q\,$ nous avons $$(q^3+q)<U^q\,,\,{\ovw}^q> +
<V^q\,,\,{\ovv}^q> + (q+q^{-1})<W^q\,,\,{\ovu}^q>\,=\,\gamma.$$

$\bullet \quad$ Enfin la $q$-sym\'etrie (\ref{pqs}) se traduit
par : $$q^2 <U^q\,,\,{\ovv}^q>=<V^q\,,\,{\ovu}^q>\,,\quad
<W^q\,,\,{\ovv}^q>=q^2<V^q\,,\,{\ovw}^q>\,,$$ $$(q^3
+q)[<U^q\,,\,{\ovw}^q> - <W^q\,,\,{\ovu}^q>]= (q^2
-1)<V^q\,,\,{\ovv}^q>.$$

Les relations de d\'ependance entre les param\`etres $\,k\,$ et
$\,\gamma\,$ d\'ecoulent de (\ref{dkg}) et (\ref{dka}) qui, en
forme explicite se traduisent par :
\begin{equation}
<K\,,\,z>=0,\quad \forall\,\,z \in \,{\osf}'^q,\quad
<z\,,\,\ovk>=0\quad \forall\,\,z \in \,{\sf V}'^q \quad
\mbox{o\`u}\label{yve}
\end{equation}
$$K\,=\,(q^3+q)uW^q +vV^q + (q+q^{-1})wU^q\,,$$
$$\ovk\,=\,(q^3+q)\,{\ovu}^q w +{\ovv}^q v + (q+q^{-1})\,{\ovw}^q
u.$$

D\'eterminons cette relation de d\'ependance en utilisant la
premi\`ere \'equation de (\ref{yve}). Il suffit de consid\'erer
$\,z={\ovu}^q\,$ dans $\,<K\,,\,z>=0$. On a alors
\begin{equation}
(q^3+q)u<W^q\,,\,{\ovu}^q> + q^2v<U^q\,,\,{\ovv}^q> + (q +
q^{-1})w<U^q\,,\,{\ovu}^q>\,= 0.\label{dkr}
\end{equation}
Or on obtient :
\begin{eqnarray}
&\displaystyle\,<U^q\,,\,{\ovu}^q>= kuu,\,<U^q\,,\,{\ovv}^q>=
kuv,\,
<V^q\,,\,{\ovw}^q>= kvw,\label{em1}\\
&\displaystyle\,<V^q\,,\,{\ovv}^q> - (q^3 + q)<W^q\,,\,{\ovu}^q>=
k(vv - (q^3 + q) wu ),\label{em2}
\end{eqnarray}
\begin{equation}
<W^q\,,\,{\ovw}^q>= kww,\,q^2 <V^q\,,\,{\ovv}^q> + (q^3 + 2q +
q{-1}) <W^q\,,\,{\ovu}^q>= \gamma.\label{em3}
\end{equation}

Les relations (\ref{em2}) et (\ref{em3})  donnent :
$$<W^q\,,\,{\ovu}^q>=\alpha[\gamma - kq^2(vv -(q^3 + q)wu)]\,,$$
$$<V^q\,,\,{\ovv}^q>=\alpha[(q^3 + q)\gamma + k(q^3 + 2q +
q^{-1})vv - k(q^3 + q)(q^3 + 2q +q^{-1})wu]$$
o\`u $\,\displaystyle\,\alpha =(2q^3 +2q+q^{-1}+q^5)^{-1}$. \\
Exprimons en outre la constante orbitale $\,c\,$ de
l'hyperbolo\"\i de quantique en fonction de $\,v\ot v\,$ et de
$\,w\ot u$. Dans l'alg\`ebre $\,\displaystyle\,{\cal
A}_{0,q}^c\,$ nous avons : $$(q^3+q)(uw - wu)+(1-q^2)vv= 0,\quad
(q^3+q)uw+vv+(q+q^{-1})wu= c \quad \mbox{donc}$$ $$q^2vv + (q^3
+2q + q^{-1})wu\,=\,c.$$ L'\'equation (\ref{dkr}) s'\'ecrit alors
:
\begin{equation}
(q^3+q)\alpha u[\gamma - kq^2(vv-(q^3+q)wu] +kq^2 vuv +
k(q+q^{-1})wuu =0.\label{ecc}
\end{equation}
Faisons  les changements  n\'ecessaires  dans  les  deux
derniers termes de l'\'equation (\ref{ecc})  de fa\c con \`a
mettre  $\,u\,$ \`a  gauche  dans  ces  termes. Nous pouvons
alors r\'e\'ecrire l'\'equation (\ref{ecc}) sous la forme :
$$\gamma\,\Gamma_1 + \Gamma_2\,vv + \Gamma_3\,wu\,=0$$ o\`u
$$\Gamma_1 = (q^3 +q)\alpha\,,$$ $$\Gamma_2 =k[-(q^3 +q)\alpha
q^2 +q^4 -q^4(q+q^{-1})\frac{q^2 -1}{q^3 +q}]\,,$$ $$\Gamma_3
=k[(q^3 +q)^2 \alpha q^2 + (q+q^{-1})].$$

Nous en d\'eduisons (pour $\,q\,$ g\'en\'erique) :
\begin{eqnarray*}
\gamma&=&-\frac{\Gamma_2}{\Gamma_1} vv -\frac{\Gamma_3}{\Gamma_1}wu \\
&=&-\frac{1}{q^2}(1+q^4)k[q^2vv + (q^3+2q+q{-1})wu] \\
&=& -\frac{1}{q^2}(1+q^4)kc,\quad c \quad \mbox{est la contante
orbitale non nulle}.
\end{eqnarray*}
On v\'erifie qu'avec cette relation de d\'ependance trouv\'ee, la
deuxi\`eme \'equation de (\ref{yve}) est satisafaite.

\subsection*{{\bf B} $\qquad$  Connexion tress\'ee}

Notons $$J_i =\frac{1}{[i]_q}Y^{i} J_{i-1},\quad i\in
\{1,\,2,\,3,\,4\}, \quad [i]_q =\frac{q^i
-q^{-i}}{q-q^{-1}},\quad \mbox{o\`u}$$ $$J_0 = \nabla_{U^q} U^q =
\alpha (uvU^q-q^2 uuV^q).$$ Par application de l'op\'erateur
$\,\displaystyle\,Y^{i}\,$ \`a la relation (\ref{nab}) nous
obtenons les quatre relations suivantes : $$-\nabla_{U^q} V^q-q^2
\nabla_{V^q} U^q\,=J_1\,,$$ $$ -\nabla_{U^q} W^q+q\nabla_{V^q}
V^q-q^4\nabla_{W^q} U^q\,=J_2\,,$$ $$\nabla_{V^q}
W^q+q^2\nabla_{W^q} V^q\,=J_3\,,$$ $$ \nabla_{W^q} W^q\,=J_4
\quad \mbox{avec :}$$ $$J_1\,=\,\alpha
\{[(q^3+q)uw-q^2vv)]U^q+(q^6+q^2-1)uvV^q-(q^3+q)uuW^q\}\,,$$
$$J_2\,=\,(1+q^2+q^4)\alpha \{-q^2vwU^q+(1-q^2)[-uw+{[
2\rbrack}^{-1}v^2]V^q + q^2uvW^q]\,$$
$$J_3\,=\,\alpha\{{-}q^3(q^2{+1})wwU^q{+}(1{+}q^4{-}q^6)vwV^q{+}q^4
[(q{+}q^{{-1}}uw{-}vv]W^q\}\,,$$
$$J_4\,=\,\alpha(wwV^q-q^4vwW^q).$$ Des relations du (\ref{qto}),
nous d\'eduisons que  : $$\nabla_{U^q} W^q\,=\,\nabla_{W^q} U^q +
\frac{q^2 -1}{q^3 +q}\nabla_{V^q} V^q +\frac{2}{q^3 +q}V^q\,,$$
$$(q^3 +q)J_2\,=\,-(q^3 +q)(1+q^4 )\nabla_{W^q} U^q+ (1+q^4
)\nabla_{V^q} V^q - 2V^q\,,$$ $$(q+q^{-1})(q^2 -1)\nabla_{W^q}
U^q+q^2 \nabla_{V^q} V^q+2V^q\,=0.$$ Par cons\'equent, nous
savons exprimer les \'el\'ements suivants  : $$\nabla_{U^q}
V^q,\,\nabla_{V^q} U^q,\,\nabla_{V^q} W^q,\,\nabla_{W^q} V^q,\,
\nabla_{W^q} U^q,\, \nabla_{U^q} W^q,\,\nabla_{V^q} V^q$$ en
fonction de $\,\alpha$. Plus pr\'ecis\'ement en fonction de
($J_1$) on a  :
\begin{eqnarray*}
\nabla_{U^q}
V^q=&\frac{1}{1+q^4}\{\{-\alpha[(q^3+q)uw-q^2vv]-2q^2\}
U^q-\alpha(q^6+q^2-1)uvV^q\\
&+\alpha(q^3+q)uuW^q\},\\
\nabla_{V^q} U^q=&\frac{1}{1+q^4}\{\{-\alpha
q^2[(q^3+q)uw-q^2vv]+2\}
U^q-\alpha q^2(q^6+q^2-1)uvV^q\\
&+\alpha q^2(q^3+q)uuW^q\} \}.
\end{eqnarray*}
En fonction de ($\,J_3\,$)  nous obtenons :
\begin{eqnarray*}
\nabla_{V^q} W^q=&\frac{1}{1+q^4}\{-\alpha q^3(1+q^2)wwU^q+
\alpha(1+q^4-q^6)vwV^q+\{\alpha q^4\\
&[(q+q^{-1})uw-vv]-2q^2\}W^q\},\\
\nabla_{W^q} V^q=&\frac{1}{1+q^4}\{-\alpha q^5(1+q^2)wwU^q+
\alpha q^2(1+q^4-q^6)vwV^q+\{\alpha q^6\\
&[(q+q^{-1})uw-vv]+2\}W^q \}.
\end{eqnarray*}
Nous avons \'egalement :
\begin{eqnarray*}
\nabla_{W^q} U^q=&\frac{1}{1+q^4}\{q^6\alpha vwU^q+\{-q^4(1-q^2)
\alpha[-uw+{[2]}^{-1}v^2]-\frac{2q}{1+q^2}\} V^q \\
&-q^6\alpha uvW^q\},\\
\nabla_{V^q} V^q=&\frac{1}{1+q^4}\{-q^3(1+q^2)^2\alpha vwU^q+
\{ q(1+q^2)(1-q^4)\alpha \{ q(1+q^2)\\
&(1-q^4)\alpha[-uw+{[2]}^{-1}v^2]+2(1-q^2)\} V^q+q^3(1+q^2)^2\alpha uvW^q\},\\
\nabla_{U^q} W^q=&\frac{1}{1+q^4}\{q^2\alpha vwU^q+\{(q^2-1)
\alpha[-uw+{[ 2]}^{-1}v^2]+\frac{2q}{1+q^2}\} V^q\\
&-q^2\alpha uvW^q\}.
\end{eqnarray*}

Il reste maintenant \`a pr\'eciser la constante
$\,\displaystyle\,\alpha$. Pour cela, consid\'erons la deuxi\`eme
\'equation du (\ref{cci})  :
\begin{equation}
[(q^3 +q)u\nabla_{W^q} +v\nabla_{V^q}
+(q+q^{-1})w\nabla_{U^q}](z)=0 \quad \forall\,z\in {\sf
V}'^q.\label{cot}
\end{equation}
Par exemple pour $\,\displaystyle\,z=U^q\,$ dans (\ref{cot}) nous
obtenons : $$\beta[(q^4 -q^2 +1)\alpha c+2][vU^q-q^2 uV^q]=0
\quad \mbox{et}\quad vU^q-q^2 uV^q\not=0 \quad \mbox{dans}\quad
T({\rm H}_q)_l.$$

D'o\`u la valeur de $\,\alpha$.

\vskip 0.5truecm

\centerline{P. Akueson} \centerline{I.S.T.V., Universit\'e de
Valenciennes} \centerline{B.P. 311 Valenciennes France}
\centerline{E-mail: akueson@univ-valenciennes.fr}


\vskip 0.5truecm

\begin{thebibliography}{99}
\bibitem[A]{A} P. Akueson : {\em \'El\'ements de G\'eom\'etrie tress\'ee}, Th\`ese
de l'Universit\'e de Valenciennes (1998).

\bibitem[AG]{AG} P. Akueson, D. Gurevich : {\em Some aspects of braided geometry :
differential calculus, tangent space, gauge theory}, to appear in
J.Phys.A., (1999).

\bibitem[BM]{BM} T. Brzezinski, S. Majid :{\em Line bundles on quantum sphere},
$q$-alg/9807052.

\bibitem[DGK]{DGK}  J.Donin, D.Gurevich, S.Khoroshkin : {\em Double quantization
of $\,{\bf C P}^n\,$ type orbits by generalized Verma modules},
JGP,28 (1998) pp.384-406.


\bibitem[DGS]{DGS} J.Donin, D.Gurevich, S.Shnider : {\em Invariant quantization
in one and two parameters on semisimple coadjoint orbits of
simple Lie groups}, J. Pure and App. Algebra 100 (1995),
pp.103-115.

\bibitem[DGR]{DGR} J.Donin, D.Gurevich, V.Rubtsov : {\em Quantum hyperboloid and
braided modules}, Alg\`ebre Non Commutative, Groupes quantiques
et invariants, Soci\'et\'e Math\'ematique de France, Collection
s\'eminaires et congr\`es, No 2 (1997), pp.103-118.

\bibitem[D]{D} J.Donin : {\em Double quantization on the coadjoint
representation of $\,sl(n)$}, Czech J.of Physics, 47, (1997), pp.
1115-1122.

\bibitem[FRT]{FRT} L.Faddeev, N.Reshetikhin, L.Takhtadhyan: {\em Quantization
of Lie groups and Lie algebras}, Leningrad Math.J.1 (1990),
pp.193-226.

\bibitem[GRR]{GRR} D.Gurevich, A. Radul, V.Rubtsov : {\em Noncommutative
differential geometry and Yang-Baxter equation}, Preprint. Publ.
Math. IHES (1991).

\bibitem[GP]{GP} D.Gurevich, D.Panyushev : {\em On Poisson pairs associated to
modified R-matrices}, Duke Math.J.73 (1994), no.1.

\bibitem[DG]{DG} D.Gurevich, J.Donin : {\em Braiding of the Lie algebra
$\,sl(2)$}, Amer.Math.Soc.Transl.(2) 167 (1995), pp.23-36.

\bibitem[GV]{GV} D.Gurevich, L.Vainerman : {\em Noncommutative analogues of
$q$-special polynomials and a $q$-integral on a sphere}, J.Phys.A.
Math.Gen.31 (1998), pp.1771-1780.

\bibitem[G1]{G1} D.Gurevich : {\em Algebraic aspects of the quantum
Yang-Baxter equation}, Leningrad. Math.J. 2 (1991), pp. 801-828.

\bibitem[G2]{G2} D.Gurevich : {\em Braided modules and reflection equations},
Quantum groups and quantum spaces. Banach Center Publ, (40),
Institut of Math, Polish Academy of Sciences, Warszawa 1997,
pp.99-109.

\bibitem[K]{K} C. Kassel : {\em Quantum groups}, Graduate texts in
mathematics, 155, (1995).

\bibitem[L]{L} J. Lambek : {\em Lectures on rings and modules}, Blaisddell
publishing company (1966).

\bibitem[LS]{LS} V. Lyubashenko, A. Sudbery : {\em Quantum Lie algebras of
type $\,A_n$}, $q$-alg /9510004.

\bibitem[M]{M} Sh. Majid : {\em Foundations of quantum groups theory},
Cambrige University
Press, (1995).

\bibitem[P]{P} P. Podl\`es : {\em Quantum spheres}, Lett.Math.Phys. 14 (1987),
pp.193-202.

\bibitem[R]{R} G. Rinehart : {\em Differential forms for general
commutative algebras},
Trans.Amer.Math.Soc.108 (1963), pp.195-222.

\bibitem[S1]{S1} J. P. Serre : {\em Modules projectifs et espaces
fibr\'es a fibre vectorielle},
exp. 23, S\'eminaire Dubreil-Pisot, Alg\`ebre et th\'eorie des
nombres, Secr\'etariat math\'ematique, Paris, 1958.

\bibitem[S2]{S2} A. Sudbery : {\em $SU_q (n)$ gauge theory}, Phys.Letters B 375
(1996), pp. 75-80.

\end{thebibliography}
\end{document}